\documentclass[12pt]{amsart}
\usepackage{amssymb}
\usepackage{amsthm}
\usepackage{amscd}

\usepackage{graphicx}
\usepackage{amsmath}
\usepackage{epic,eepic}
\theoremstyle{plain}
\newtheorem{lemma}{Lemma}[section]
\newtheorem{prop}[lemma]{Proposition}
\newtheorem{thm}[lemma]{Theorem}
\newtheorem{cor}[lemma]{Corollary}
\newtheorem{aplemma}{Lemma~A.\hspace{-1.5mm}}
\newtheorem{approp}{Proposition~A.\hspace{-1.5mm}}
\newtheorem{apthm}{Theorem~A.\hspace{-1.5mm}}
\newtheorem{apcor}{Corollary~A.\hspace{-1.5mm}}

\usepackage[all]{xy}

\newtheorem{conj}[lemma]{Conjecture}

\theoremstyle{definition}

\newtheorem{rema}[lemma]{Remark}

\newtheorem{remb}{Remark}

\newtheorem{defi}[lemma]{Definition}
\newtheorem{exa}[lemma]{Example}
\newtheorem{aprem}{Remark~A.\hspace{-1.5mm}}
\newtheorem{apdefi}{Definition~A.\hspace{-1.5mm}}
\newcommand{\bde}{\begin{defi}}
\newcommand{\ede}{\end{defi}\vspace{1mm}}
\newcommand{\ble}{\begin{lemma}}
\newcommand{\ele}{\end{lemma}}
\newcommand{\bpr}{\begin{prop}}
\newcommand{\epr}{\end{prop}}
\newcommand{\bt}{\begin{thm}}
\newcommand{\et}{\end{thm}}
\newcommand{\bco}{\begin{cor}}
\newcommand{\eco}{\end{cor}}
\newcommand{\bre}{\begin{rema}}
\newcommand{\ere}{\end{rema}}
\newcommand{\brea}{\begin{rema}}
\newcommand{\erea}{\end{rema}\vspace{1mm}}
\newcommand{\breb}{\begin{remb}}
\newcommand{\ereb}{\end{remb}\vspace{1mm}}
\newcommand{\bex}{\begin{exa}}
\newcommand{\eex}{\end{exa}}
\newcommand{\bpf}{\begin{proof}}
\newcommand{\epf}{\end{proof}\vspace{1mm}}

\newcommand{\bade}{\begin{apdefi}}
\newcommand{\eade}{\end{apdefi}}
\newcommand{\bale}{\begin{aplemma}}
\newcommand{\eale}{\end{aplemma}}
\newcommand{\bapr}{\begin{approp}}
\newcommand{\eapr}{\end{approp}}
\newcommand{\bat}{\begin{apthm}}
\newcommand{\eat}{\end{apthm}}
\newcommand{\baco}{\begin{apcor}}
\newcommand{\eaco}{\end{apcor}}
\newcommand{\bare}{\begin{aprem}}
\newcommand{\eare}{\end{aprem}}

\newcommand{\be}{\begin{enumerate}}
\newcommand{\ee}{\end{enumerate}}
\newcommand{\bcd}{\[\begin{CD}}
\newcommand{\ecd}{\end{CD}\]}
\newcommand{\bit}{\begin{itemize}}
\newcommand{\eit}{\end{itemize}}
\newcommand{\bq}{\begin{quote}}
\newcommand{\eq}{\end{quote}}
\newcommand{\ba}{\begin{array}}
\newcommand{\ea}{\end{array}}

\newcommand{\mcD}{\mathcal{D}}
\newcommand{\mcE}{\mathcal{E}}
\newcommand{\mcF}{\mathcal{F}}

\newcommand{\mcH}{\mathcal{H}}
\newcommand{\mcI}{\mathcal{I}}

\newcommand{\mcL}{\mathcal{L}}

\newcommand{\mcN}{\mathcal{N}}
\newcommand{\mcO}{\mathcal{O}}
\newcommand{\mcP}{\mathcal{P}}

\newcommand{\mcS}{\mathcal{S}}
\newcommand{\mcT}{\mathcal{T}}

\newcommand{\mcW}{\mathcal{W}}

\newcommand{\mbA}{\mathbb{A}}

\newcommand{\mbC}{\mathbb{C}}

\newcommand{\mbG}{\mathbb{G}}

\newcommand{\mbL}{\mathbb{L}}

\newcommand{\mbP}{\mathbb{P}}

\newcommand{\mbZ}{\mathbb{Z}}

\newcommand{\mfB}{\mathfrak{B}}
\newcommand{\mfC}{\mathfrak{C}}

\newcommand{\mfI}{\mathfrak{I}}

\newcommand{\mfO}{\mathfrak{O}}
\newcommand{\mfP}{\mathfrak{P}}
\newcommand{\mfQ}{\mathfrak{Q}}

\newcommand{\mfS}{\mathfrak{S}}

\newcommand{\mfp}{\mathfrak{p}}

\newcommand{\migi}{\rightarrow}
\newcommand{\longmigi}{\longrightarrow}

\newcommand{\isom}{\stackrel{\sim}{\migi}}

\newcommand{\migiincl}{\hookrightarrow}

\newcommand{\migisurj}{\twoheadrightarrow}

\newcommand{\A}{\mbA}

\newcommand{\mr}{\mathrm}
\newcommand{\hidden}[1]{\,}

\begin{document}

\title[Quantization on Algebraic Curves  with Frobenius-Projective Structure]{Quantization on Algebraic Curves \\ with Frobenius-Projective Structure}
\author{Yasuhiro Wakabayashi}
\date{}
\markboth{Yasuhiro Wakabayashi}{}
\maketitle
\footnotetext{Y. Wakabayashi: Department of Mathematics, Tokyo Institute of Technology, 2-12-1 Ookayama, Meguro-ku, Tokyo 152-8551, JAPAN;}
\footnotetext{e-mail: {\tt wkbysh@math.titech.ac.jp};}
\begin{abstract}
In the present paper, 
we study 
the relationship between deformation quantizations and Frobenius-projective structures defined on an algebraic curve in positive characteristic.
A Frobenius-projective structure is an analogue  of a complex  projective structure on a Riemann surface, which was introduced by Y. Hoshi.
Such an additional structure has some equivalent objects, e.g., a dormant $\mr{PGL}_2$-oper and a projective connection having a full set of solutions.
The main result of the present paper provides a canonical construction of a Frobenius-constant quantization on the cotangent space minus the zero section on an algebraic curve  by means of a Frobenius-projective structure.
It may  be thought of as a positive characteristic analogue of a result by D. Ben-Zvi and I. Biswas.
Finally, we give a higher-dimensional variant of this result, as proved by I. Biswas in the complex case.

\end{abstract}
\tableofcontents 

\section*{Introduction}

A (deformation) quantization of a symplectic manifold is a noncommutative deformation of the structure sheaf which is, in a certain sense, compatible with
the symplectic structure.
We know that every symplectic manifold admits quantizations (cf. ~\cite{DWL1}; ~\cite{DWL2}; ~\cite{Fe}), but in general 
the quantization is neither unique nor canonically constructed. 
Therefore, construction of quantizations will be one of the important subjects in this theory.
In  \cite{BZBi},
 D. Ben-Zvi and I. Biswas provided a canonical construction
 of a quantization on (the total space of) the cotangent bundle minus the zero section on a given Riemann surface $C$ equipped with a projective structure.
Here, recall that a projective structure on $C$ is an atlas of  coordinate charts defining $C$  whose transition functions may be expressed as Mobi\"{u}s transformations.
Projective structures 
have a very major role to play in understanding the framework of uniformization theorem of Riemann surfaces 
 and have mutually equivalent objects, including $\mr{PGL}_2$-opers and projective connections, etc..
Every Riemann surface $C$  admits a projective structure, and the space of all projective structures on $C$ forms an affine space for the space $H^0 (C, \Omega_{C}^{\otimes 2})$ of quadratic differentials on $C$.
The idea behind the construction of D. Ben-Zvi and I. Biswas
is that the canonical construction of a quantization on the complex projective line $\mbP^1_\mbC$,  which is invariant under the action of $\mr{PGL}_2(\mbC)$ (= the group of Mobi\"{u}s transformations on $\mbP^1_\mbC$), may extend naturally to any Riemann surface once we choose a projective structure.

On algebraic curves in positive characteristic, there are analogous objects of  projective structures, called {\it Frobenius-projective structures}. 
 The notion of a Frobenius-projective structure was
  introduced by Y. Hoshi (cf. ~\cite{Hos}, \S\,2, Definition 2.1) as a certain collection of locally defined \'{e}tale maps on a prescribed curve to the projective line.
Just as in the complex case, any  smooth curve in positive characteristic  admits such a structure.
Also, Frobenius-projective structures are equivalent to, e.g., {\it dormant}  $\mr{PGL}_2$-opers and projective connections {\it having a full set of solutions}.
That is to say, given a connected smooth curve $X$ in characteristic $p>2$ and a theta characteristic $\mbL := (\mcL, \psi_\mcL : \mcL^{\otimes 2} \isom \Omega_X)$ (cf. \S\,\ref{S8}), we obtain
   the following diagram consisting of bijective correspondences in parallel with 
the classical result on Riemann surfaces: 
\begin{align}
\vcenter{\xymatrix@C=16pt@R=16pt{
 & \mfP \mfS_X^\mr{F} \ar[ldd] \ar[rdd]&
 \\
 & & \\
 \mfO \mfp^{^\mr{Zzz...}}_{\mr{PGL}_2, X} \ar[rr]_{\sim} \ar[uur]^{\rotatebox{40}{{\tiny $\sim$}}}& & \mfP \mfC^{2, \mr{full}}_{X, \mbL}, \ar[uul]_{\rotatebox{-40}{{\tiny $\sim$}}} \ar[ll]
}} 
\end{align}
where
\begin{itemize}
\item[] \hspace{-3mm} $\mfP \mfS_X^\mr{F}$ := the set of Frobenius-projective structures on $X$ (cf. (\ref{Efllo2})); 
\vspace{1mm}
\item[] \hspace{-3mm} $\mfO \mfp^{^\mr{Zzz...}}_{\mr{PGL}_2, X}$ := the set of isomorphism classes of dormant $\mr{PGL}_2$-opers on $X$ (cf. (\ref{E69207}));
\vspace{1mm}
\item[] \hspace{-3mm} $\mfP \mfC^{2, \mr{full}}_{X, \mbL}$ := the set of projective connections for $\mbL$
having a full set of solutions (cf. (\ref{E41180})).
\end{itemize}
(We also discuss, in the present paper, certain intermidiate   objects equivalent to them, called {\it dormant $(\mr{SL}_2, \mbL)$-opers}.)

The purpose of the present paper is to prove an analogous assertion of  D. Ben-Zvi and I. Biswas, i.e., a canonical construction of a quantization by means of  the curve $X$ together with a choice among such additional structures.
In ~\cite{BK1}, ~\cite{Kon}, and  ~\cite{Y}, it has been shown that the  theory of quantizations can be made to work in the algebraic setting.
Also, we can find, in ~\cite{BeKa1} (and  ~\cite{BeKa2}), the study of a special class of quantizations on symplectic algebraic varieties in positive characteristic, called {\it Frobenius-constant quantizations}.
They are quantizations with large center  in some suitable sense, and has a cohomological classification given in the point of view of formal geometry.

Let $\mbA (\Omega_X)^\times$ denote  the complement of the zero section  in (the total space of) the cotangent bundle of $X$; it admits
a symplectic structure $\check{\omega}^\mr{can}$ defined as one half of the Liouville symplectic form.
Thus, it makes sense to speak of a (Frobenius-constant) quantization on the symplectic variety  $(\mbA (\Omega_X)^\times, \check{\omega}^\mr{can})$.
Denote (cf. (\ref{Er33})) by
\begin{align}
\mfQ^{\mr{FC}}_{(\mbA (\Omega_X)^\times, \check{\omega}^\mr{can})}
\end{align}
the set of Frobenius-constant quantizations on  
$(\mbA (\Omega_X)^\times, \check{\omega}^\mr{can})$.
Then, the main result of the present paper (cf. Theorem \ref{T0090}) provides a  canonical {\it injective} assignment from
a Frobenius-projective structure (or equivalently, a dormant indigenous bundle, or a projective connection having a full set of solutions) to a Frobenius-constant quantization on $(\mbA (\Omega_X), \check{\omega}^\mr{can})$,  as displayed below: 
\begin{align}
\hspace{-5mm}\vcenter{\xymatrix@C=16pt@R=16pt{
 & \mfP \mfS_X^\mr{F} \ar[ldd] \ar[rdd]&
 \\
 & & \\
 \mfO \mfp^{^\mr{Zzz...}}_{\mr{PGL}_2, X} \ar[rr]_{\sim} \ar[uur]^{\rotatebox{40}{{\tiny $\sim$}}}& & \mfP \mfC^{2, \mr{full}}_{X, \mbL}, \ar[uul]_{\rotatebox{-40}{{\tiny $\sim$}}} \ar[ll]
}} 
\hspace{15mm} \scalebox{2}{$\rightsquigarrow$} \hspace{15mm} \mfQ^\mr{FC}_{(\mbA(\Omega_X)^\times, \check{\omega}^\mr{can})}.
\end{align}
In particular, we can think of $\mfP \mfS^\mr{F}_X$ as a subset of $\mfQ^\mr{FC}_{(\mbA(\Omega_X)^\times, \check{\omega}^\mr{can})}$ via this assignment,  and hence, give a lower bound of the number of Frobenius-constant quantizations on $(\mbA (\Omega_X)^\times, \check{\omega}^\mr{can})$ by applying the result in ~\cite{Wak}.

The present paper is organized as follows.
The first section contains the necessary definitions
 and conventions used in our discussion, including a symplectic structure, a differential operator, and a theta characteristic.
In the second section, we recall the notion of a Frobenius-constant quantization and discuss some related topics.
For instance, it is observed that Frobenius-constant quantizations are functorial  with respect to pull-back via \'{e}tale morphisms (cf. \S\,\ref{S10}) and have a descent property via  finite Galois coverings (cf. \S\,\ref{S41}).
In the third section, we discuss various bijective correspondences between
 Frobenius-projective structures on a curve  and some equivalent objects, i.e., dormant $\mr{PGL}_2$-opers, dormant $(\mr{SL}_2, \mbL)$-opers, and  projective connections with a full set of solutions.
Some results mentioned in that section have been  essentially obtained in other literatures, e.g.,  ~\cite{Hos} (which gives $\mfP \mfS^\mr{F}_X \isom \mfO \mfp^{^\mr{Zzz...}}_{\mr{PGL}_2, X}$) and ~\cite{BD2} together with ~\cite{Kat2}
 (which gives $\mfO \mfp^{^\mr{Zzz...}}_{(\mr{SL}_n, \mbL), X}\isom \mfP \mfC^{n, \mr{full}}_{X, \mbL}$).
Also, in ~\cite{Wak2}, we can find generalizations of these correspondences to a family of pointed stable curves. 
But, unfortunately, many of them seem not to be standard and are unavoidable when we complete the proof of the main theorem, so we decided to review them here  and the contents became nearly self-contained.
The fourth section is devoted to state and prove the main theorem.
As carried out in ~\cite{BZBi},  we first construct, by means of a Frobenius projective structure, a Frobenius-constant quantization on the complement of the zero section in the total space of $\mcL$.
This quantization turns out to be invariant under the natural involution, and hence, descends to a Frobenius-constant quantization on $(\mbA (\Omega_X)^\times, \check{\omega}^\mr{can})$.
Moreover, the injectivity of  this assignment is proved by examining  the behavior of the noncommutative multiplication in each quantization. 
In the final section, 
we discuss (cf. Theorem \ref{T009h1gh}) a higher-dimensional variant of our main theorem, which may be thought of as a positive characteristic analogue of a result in ~\cite{Bi}.

\vspace{5mm}
\hspace{-4mm}{\bf Acknowledgement} \leavevmode\\
 \ \ \ 
The author would like to thank all algebraic curves  equipped with a Frobenius-projective structure, who live in the world of mathematics, for their useful comments and heartfelt encouragement. The author was partially supported by the Grant-in-Aid for Scientific Research (KAKENHI No. 18K13385).

\vspace{10mm}
\section{Preliminaries} \vspace{3mm}

In this section, we prepare the notation and conventions used in the present paper.
Throughout the present paper, let us fix an odd prime $p$ and   an algebraically closed field $k$ of characteristic $p$.
Unless otherwise stated, all schemes and morphisms of schemes are implicitly assumed to be over $k$, and products of schemes are taken over $k$.
We use the word {\it  variety} (resp., {\it curve}) to mean a  finite type integral scheme over $k$ (resp., a  finite type integral scheme over $k$ of dimension $1$).
For each positive integer $n$, we shall write
$\mbA^n$ (resp., $\mbP^n$) for the affine space (resp., the projective space) over $k$ of dimension $n$.
Also,  write
$\mbA^{n \times} := \mbA^n \setminus \{0\}$.

\vspace{5mm}
\subsection{Vector bundles} \label{S7}
\leavevmode\\ \vspace{-4mm}

Let 
$S$ be a smooth variety of dimension $n \geq 0$.  
Given a vector bundle $\mcF$ on $S$ (i.e., a locally free coherent sheaf on $S$),
we denote by $\mcF^\vee$ its dual sheaf, i.e., $\mcF^\vee := \mcH om_{\mcO_X} (\mcF, \mcO_X)$.
Let  $\mbA (\mcF)$ and $\mbP(\mcF)$ denote the relative affine  and projective spaces respectively associated with $\mcF$, i.e.,
\begin{align} \label{E04}
\A (\mcF) := \mcS pec (S^\bullet (\mcF^\vee)), \ \ \ \mbP (\mcF) := \mcP roj (S^\bullet (\mcF^\vee)),
\end{align}
 where $S^\bullet (\mcF^\vee)  \ \left(:= \bigoplus_{i \geq 0}S^i (\mcF^\vee)\right)$ denotes the symmetric algebra over $\mcO_S$ associated with  $\mcF^\vee$.
 Also, write
 \begin{align}
 \A (\mcF)^\times 
 \end{align}
for  the complement of the zero section $S \migi \mbA (\mcF)$ in $\A (\mcF)$, which admits a natural projection
\begin{align} \label{E06}
\pi_\mcF :   \A (\mcF)^\times  \migi  \mbP (\mcF) 
\end{align}
over $S$.

We shall write $\Omega_{S}$  for  the sheaf of $1$-forms (in other words, the cotangent bundle) on 
$S$ relative to $k$ and $\mcT_S$ for its dual.
By the smoothness assumption on $S$, both $\Omega_S$ and $\mcT_S$ turn out to be  vector bundles
 of rank $n$.
We write
$d : \mcO_S \migi \Omega_S$
for the universal derivation.
Moreover, denote by $\omega_{S}$ the canonical line bundle of $S$ (relative to $k$), which is canonically isomorphic to the determinant line bundle $\mr{det}(\Omega_S) := \bigwedge^n \Omega_S$ of $\Omega_S$.

\vspace{5mm}
\subsection{Symplectic structures} \label{S700}
\leavevmode\\ \vspace{-4mm}

Recall that a {\bf symplectic structure}  on $S$  is a nondegenerate closed $2$-form $\omega \in \Gamma (S, \bigwedge^2\Omega_{S})$.
Here, we  say that $\omega$ is {\it nondegenerate}  if the morphism $\Omega_{S} \migi \mcT_{S} \ \left(= \Omega_{S}^\vee \right)$ induced naturally  by $\omega$ is an isomorphism.
A {\bf symplectic variety (over $k$)} is  a pair $(S, \omega)$ consisting of a   smooth variety $S$   and  a symplectic structure $\omega$ on it.
An {\bf isomorphism} $(S, \omega) \isom (S', \omega')$ between symplectic varieties is an isomorphism $S \isom S'$ preserving the respective symplectic structures.

As is well-known,  the variety $\A(\Omega_S)$ (i.e., the total space of the cotangent bundle of $S$)   has  a canonical symplectic structure
\begin{align} \label{W104}
\omega^\mr{can}_S\in \Gamma (\A (\Omega_S), {\bigwedge}^2 \Omega_{\A (\Omega_S)}).
\end{align}
often called the {\it Liouville symplectic form}.
If there is no fear of causing confusion, we write $\omega^\mr{can}$ instead of $\omega^\mr{can}_S$ for simplicity.
If $q_1, \cdots, q_{n}$ are local coordinates in $S$ and  $q^\vee_1, \cdots,  q^\vee_{n}$  denote  the dual coordinates in $\A (\Omega_{S})$, then  $\omega^{\mr{can}}$ may be expressed locally as
$\omega^{\mr{can}} = \sum_{i=1}^{n} dq^\vee_i \wedge dq_i$.
By abuse of notation, we also use the notation  $\omega^\mr{can}$ to denote  the restriction of $\omega^\mr{can}$ to the open subscheme $\A(\Omega_S)^\times \ \left(\subseteq \A (\Omega_S) \right)$. 
Also, for each $c \in k^\times$, $c \cdot \omega^\mr{can}$ forms a symplectic structure. 
In particular, by letting $\check{\omega}^\mr{can} := \frac{1}{2}\cdot \omega^\mr{can}$,
 we have symplectic varieties
\begin{align} \label{Errgt8}
(\A(\Omega_S), \check{\omega}^\mr{can}), \hspace{5mm}
 (\A(\Omega_S)^\times, \check{\omega}^\mr{can}).
\end{align}

\vspace{5mm}
\subsection{Differential operators} \label{S83}
\leavevmode\\ \vspace{-4mm}

 We shall recall the notion of a differential operator.
Let $\mcL_i$ ($i=1,2$) be line bundles on $S$.
By a {\it differential operator} from $\mcL_1$ to $\mcL_2$, we mean a $k$-linear  morphism $D  :\mcL_1 \migi \mcL_2$ 
locally expressed, after fixing identifications $\mcL_1 \cong \mcL_2 \cong \mcO_S$ and a local coordinate system $\vec{x}:= (x_1, \cdots, x_n)$ in $S$, as 
\begin{align} \label{E442}
D 
: v \mapsto 
D(v) =  \sum_{\alpha \in \mbZ_{\geq 0}^n} a_\alpha \cdot \partial^\alpha_{\vec{x}} (v)
\end{align}
 by means of some local sections $a_\alpha \in \mcO_S$  with $a_\alpha =0$ for almost all $\alpha$,
  where 
 for each $\alpha := (\alpha_1, \cdots, \alpha_n) \in \mbZ_{\geq 0}^n$, 
  we write
 $\partial^\alpha_{\vec{x}}(v) := \frac{\partial^{|\alpha|}(v)}{\partial x_1^{\alpha_1} \cdots \partial x_n^{\alpha_n}}$ ($|\alpha| := \alpha_1 + \cdots + \alpha_n$).
If $a_\alpha=0$ for any $\alpha$ with $|\alpha| \geq p$,
then   $j_\mr{max} := \mr{max}\left\{ |\alpha| \, | \, a_\alpha \neq 0 \right\} \ (<p)$ is well-defined (i.e., depend only on $D$,  not on the choice of the local expression (\ref{E442})).
In this situation, 
 we say that $D$ is {\it of order $j_\mr{max}$}. (We say that $D$ is {\it of order $-\infty$} if $D=0$.)

Given a nonnegative integer $j$ with $j<p$, we denote by
\begin{align}
\mcD \textit{iff}^{\,\leq j}_{\mcL_1, \mcL_2}
\end{align} 
the Zariski sheaf on $S$ consisting of locally defined differential operators from $\mcL_1$ to $\mcL_2$ of order $\leq j$; it is a subsheaf of the sheaf $\mcH om_k (\mcL_1, \mcL_2)$ of locally  defined $k$-linear morphisms $\mcL_1 \migi \mcL_2$.
In the case where $\mcL_1 = \mcL_2 =\mcO_X$, we write
\begin{align}
\mcD_X^{\leq i} := \mcD \textit{iff}^{\,\leq j}_{\mcO_X, \mcO_X}.
\end{align}
Note that $\mcD \textit{iff}^{\,\leq j}_{\mcL_1, \mcL_2}$
admits two different structures of $\mcO_S$-module --- one as given by left multiplication (where we denote this $\mcO_S$-module by 
${^l \mcD} \textit{iff}^{\,\leq j}_{\mcL_1, \mcL_2}$),
and the other given by right multiplication (where we denote this $\mcO_S$-module by ${^r \mcD} \textit{iff}^{\,\leq j}_{\mcL_1, \mcL_2}$) ---.
Given an $\mcO_X$-module $\mcF$, we equip the tensor product $\mcF \otimes {^l \mcD}_{X}^{\leq j}$ (resp., ${^r \mcD}_X^{\leq j} \otimes \mcF$) with an $\mcO_X$-module  structure  arising from 
  ${^r \mcD}_X^{\leq j}$ (resp., ${^l \mcD}_X^{\leq j}$).
Then,
the composition with the $k$-linear morphism $\mcL_2 \otimes \mcD_X^{\leq j} \migi \mcL_2$ given by
$v \otimes D \mapsto v \otimes D(1)$ yields
  an identification
\begin{align} \label{Edr56}
\mcH om_{\mcO_X} (\mcL_1, \mcL_2 \otimes \mcD_X^{\leq j}) \isom \mcD {\it iff}_{\mcL_1, \mcL_2}^{\leq j}.
\end{align}

Moreover, the assignment $D =  \sum_{\alpha \in \mbZ_{\geq 0}^n} a_\alpha \cdot \partial^\alpha_{\vec{x}} \mapsto  \sum_{|\alpha| =j} a_\alpha \cdot \partial^\alpha_{\vec{x}}$ gives a well-defined isomorphism  of $\mcO_S$-modules
\begin{align}
{ \mcD} \textit{iff}^{\,\leq j}_{\mcL_1, \mcL_2}/{\mcD} \textit{iff}^{\,\leq (j-1)}_{\mcL_1, \mcL_2}
\isom \mcH om_{\mcO_S} (\mcL_1, \mcL_2 \otimes S^j (\mcT_{S})),
\end{align}
where $S^j(\mcT_{S})$ denotes the $j$-th component of the symmetric power of $\mcT_S$.
Denote by $\Sigma$ the composite
\begin{align}
\Sigma : \mcD \textit{iff}^{\,\leq j}_{\mcL_1, \mcL_2} \migisurj
\mcD \textit{iff}^{\,\leq j}_{\mcL_1, \mcL_2}/\mcD \textit{iff}^{\,\leq (j-1)}_{\mcL_1, \mcL_2} \isom \mcH om_{\mcO_S} (\mcL_1, \mcL_2 \otimes S^j (\mcT_{S})).
\end{align}
For each local section $D \in \mcD \textit{iff}^{\,\leq j}_{\mcL_1, \mcL_2}$, we refer to $\Sigma (D)$ as the {\it principal symbol} of $D$.

Next,  let us write 
\begin{align}
S^{(1)}
\end{align}
 for  the Frobenius twist of $S$ over $k$  (i.e., the base-change of $S$ via the absolute Frobenius morphism of $k$) 
   and 
   \begin{align}
   F_{S/k}: S \migi S^{(1)}
   \end{align}
    for the relative Frobenius morphism of $S$ over $k$.
To simplify the notation,
we regard each $\mcO_{S^{(1)}}$-module  (resp., $\mcO_S$-module)
as a sheaf on $S$ (resp., on $S^{(1)}$) via the underlying homeomorphism of $F_{S/k}$.

Notice that each differential operator  $D : \mcL_1 \migi \mcL_2$ of order $j$ may be
 considered as an {\it $\mcO_{S^{(1)}}$-linear} morphism 
 $F_{S/k*}(\mcL_1) \migi F_{S/k*}(\mcL_2)$ 
  via the underlying homeomorphism of $F_{S/k}$.
It follows that the kernel $\mr{Ker}(D)$ forms an $\mcO_{S^{(1)}}$-submodule of $F_{S/k*}(\mcL_1)$.

\vspace{3mm}
\bde \label{D01fgh} \leavevmode\\
 \ \ \ 
We shall say that {\bf $D$ has a full set of  solutions} if  $\mr{Ker}(D)$  is a vector bundle on $S^{(1)}$ of  rank $j$.
\ede

\vspace{5mm}
\subsection{Theta characteristics} \label{S8}
\leavevmode\\ \vspace{-4mm}

 By a {\bf theta characteristic} on $S$, we mean
a pair 
\begin{align}
\mbL := (\mcL, \psi_\mcL)
\end{align}
 consisting of a line bundle $\mcL$ on $S$  and an isomorphism $\psi_\mcL : \mcL^{\otimes (n+1)} \isom \omega_S$ between line bundles.
As is well-known, 
any smooth curve always admits a theta characteristic.

\vspace{3mm}
\begin{exa} \label{E5d10}
\leavevmode\\
\ \ \
We shall observe that there exists a canonical theta characteristic on the projective space $\mbP^{n} = \mr{Proj}(k[x_0, x_1, \cdots, x_n])$.
Let
\begin{align} \label{Effgh}
\eta_0 : \mcO_{\mbP^n}(-1) \migi \mcO_{\mbP^n}^{\oplus (n+1)}
\end{align}
be the $\mcO_{\mbP^n}$-linear injection given by $w \mapsto \sum_{i=0}^n w x_i\cdot  e_i$ for each local section $w \in \mcO_{\mbP^n}(-1)$, where $(e_0, \cdots, e_n)$ is a canonical basis of $\mcO_{\mbP^n}^{\oplus (n+1)}$.
The composite
\begin{align}
\mcO_{\mbP^n}(-1) \xrightarrow{\eta_0}  \mcO_{\mbP^n}^{\oplus (n+1)} \xrightarrow{d^{\oplus (n+1)}}  \Omega_{\mbP^n} \otimes \mcO_{\mbP^n}^{\oplus (n+1)} \migisurj \Omega_{\mbP^n} \otimes (\mcO_{\mbP^n}^{\oplus (n+1)}/\mr{Im}(\eta_0)), 
\end{align}
which is verified to be $\mcO_{\mbP^n}$-linear, induces an isomorphism of $\mcO_{\mbP^n}$-modules
\begin{align} \label{E5568}
\mcO_{\mbP^n}(-1) \otimes (\mcO_{\mbP^n}^{\oplus (n+1)}/\mr{Im}(\eta_0))^\vee \isom \Omega_{\mbP^n}. 
\end{align}
Moreover,
we have a composite isomorphism
\begin{align} \label{G3345}
\mcO_{\mbP^n} \isom \mr{det} (\mcO_{\mbP^n}^{\oplus (n+1)}) 
 \isom \mcO_{\mbP^n}(-1)\otimes \mr{det} (\mcO_{\mbP^n}^{\oplus (n+1)}/\mr{Im}(\eta_0)),
\end{align}
where the first isomorphism is given by $1 \mapsto e_0 \wedge \cdots \wedge e_n$ and the second isomorphism arises from 
 the short exact sequence 
\begin{align}
0 \xrightarrow{}  \mcO_{\mbP^n}(-1) \xrightarrow{\eta_0} \mcO_{\mbP^n}^{\oplus (n+1)} \xrightarrow{} \mcO_{\mbP^n}^{\oplus (n+1)}/\mr{Im}(\eta_0) \xrightarrow{} 0.
\end{align}
Denote by $\psi_0$ the composite isomorphism
\begin{align}
\psi_0 : \mcO_{\mbP^n}(-1)^{\otimes (n+1)} \ &\left(=\mcO_{\mbP^n}(-n) \otimes \mcO_{\mbP^n}(-1) \right) \\
& \, \isom  \mcO_{\mbP^n}(-n) \otimes \mr{det}(\mcO_{\mbP^n}^{\oplus (n+1)}/\mr{Im}(\eta_0))^\vee \notag \\
& \,  \isom \mr{det}(\mcO_{\mbP^n}(-1) \otimes (\mcO_{\mbP^n}^{\oplus (n+1)}/\mr{Im}(\eta_0))^\vee) \notag \\
& \,  \isom \left(\mr{det}(\Omega_{\mbP^n}) = \right) \ \omega_S, \notag
\end{align}
where the first isomorphism follows from (\ref{G3345}) and the third isomorphism follows from (\ref{E5568}).
Thus, we have obtained a theta characteristic
\begin{align} \label{E09076}
\mbL_0 := (\mcO_{\mbP^n}(-1), \psi_0)
\end{align}
on $\mbP^n$.
\end{exa}
\vspace{3mm}

\vspace{10mm}
\section{Frobenius-constant quantizations} \vspace{3mm}

In this section, we recall the notion of a Frobenius-constant (= FC) quantization on a given symplectic variety and discuss some related topics.

\vspace{5mm}
\subsection{Quentizations} \label{S4}
\leavevmode\\ \vspace{-4mm}

Let $(S, \omega_S)$ be a symplectic variety.
The  nondegenerate pairing $\mcT_{S} \otimes_{\mcO_S} \mcT_{S} \migi \mcO_S$ given by $\omega_S$
becomes a pairing
$\omega^{-1}_S : \Omega_{S} \otimes_{\mcO_S} \Omega_{S} \migi \mcO_S$ via  $\Omega_S \isom \mcT_S$ induced by $\omega_S$.
Thus, we obtain  a skew-symmetric  $k$-bilinear map
\begin{align} 
\{-, - \}_{\omega} : \mcO_S \times \mcO_S \migi \mcO_S
\notag \end{align}
 defined by 
$\{ f, g \}_{\omega} := \omega^{-1}_S(df, dg)$.
One verifies from the closedness of  $\omega_S$ that $\{-, -\}_\omega$ defines  a Poisson bracket in the usual sense.
Here, let $k[[\hslash]]$ denote the ring of formal power series in the  variable $\hslash$ over $k$, and write $\mcO_S [[\hslash ]] := \varprojlim_{j \geq 1} \mcO_S [\hslash]/(\hslash^j)$.
In this article, we shall define 
a {\bf quantization} on   $(S, \omega_S)$ to be 
a sheaf of (noncommutative) flat $k[[\hslash ]]$-algebras $\mcO^\hslash_S$ on $S$ 
such that $\mcO^\hslash_S = \mcO_S [[\hslash ]]$ (as an equality of sheaves of $k[[\hslash]]$-modules)   
 and
the
 commutator in $\mcO^\hslash_S$ is equal to $\hslash \cdot \{-,-\}_\omega$ mod $\hslash^2 \cdot \mcO_S^\hslash$.
Moreover,
a {\bf Frobenius-constant quantization} (or,  an {\bf FC quantization},  for short) on $(S, \omega_S)$  (cf.  ~\cite{BeKa1}, Definition 3.3; \cite{BeKa2}, Definitions 1.1 and 1.4) is
a quantizaton $\mcO^\hslash_S$ on $(S, \omega_S)$ such that 
the image of the natural inclusion 
$\mcO_{S^{(1)}} [[\hslash]]\migiincl \mcO_S[[\hslash ]] \ \left(= \mcO^\hslash_S\right)$ coincides with
the center $Z(\mcO^\hslash_S)$ of $\mcO_S^\hslash$.
We shall write
\begin{align} \label{Er33}
\mfQ_{(S, \omega_S)}^{\mr{FC}}
\end{align}
for the set  of FC quantizations on $(S, \omega_S)$.

\vspace{5mm}
\subsection{Pull-back of FC quantizations} \label{S10}
\leavevmode\\ \vspace{-4mm}

If we are given an FC quantization on a prescribed symplectic variety,
then it induces an FC quantization on each open subvariety via restriction.
More generally, 
 we can construct  the pull-back of an FC quantization via an {\it \'{e}tale} morphism, as follows.
 
 Let
   $(T, \omega_T)$ be another symplectic variety
    and $f : T \migi S$ an \'{e}tale morphism 
   with $f^*(\omega_S) = \omega_T$.
The \'{e}taleness of $f$ implies that
 the commutative square diagram
 \begin{align} \label{E073}
 \xymatrix{
 T \ar[r]^{f} \ar[d]_{F_{T/k}} & \ar[d]^{F_{S/k}} S
 \\
 T^{(1)} \ar[r]_{f^{(1)}} &  S^{(1)}
 }
\end{align}
is cartesian, where $f^{(1)}$ denotes  the base-change of $f$  via the absolute Frobenius morphism of $k$.
Given an FC quantization $\mcO_{S}^\hslash$
 on $(S, \omega_S)$,
we set
\begin{align}
f^*  (\mcO_{S}^\hslash) := \varprojlim_{j >0} \left(\mcO_{T^{(1)}}\otimes_{f^{-1}(\mcO_{S^{(1)}})} f^{-1}(\mcO_{S}^\hslash/(\hslash^j) \right).
\end{align}
Then, since
$Z (f^*  (\mcO_{S}^\hslash)) \cong \mcO_{T^{(1)}} \otimes_{f^{-1}(\mcO_{S^{(1)}})} \varprojlim_{j>0}\left(Z (\mcO_S^\hslash)/(\hslash^j) \right)$,
the sheaf  $f^*  (\mcO_{S}^\hslash)$ specifies an FC quantization on $(T, \omega_T)$.
We shall refer to $f^*  (\mcO_{S}^\hslash)$ as the {\bf pull-back} of $\mcO_S^\hslash$ (via $f$).

\vspace{5mm}
\subsection{Equivariant FC quantizations} \label{S41}
\leavevmode\\ \vspace{-4mm}

Let $(S, \omega_S)$ be as above
 and $G$ a finite group acting freely on $(S, \omega_S)$, i.e., $S$ is equipped with a free $G$-action preserving $\omega_S$.  
 
 A {\bf Frobenius $G$-constant quantization} (or, a {\bf $G$-FC quantization}) on $(S, \omega_{S})$ (cf.  ~\cite{BeKa1}, Definition 5.5) is an FC quantization
 $\mcO_S^\hslash$ on $(S, \omega_S)$ compatible, in the natural sense, with the $G$-action on $S$.
Denote by 
\begin{align}
\mfQ^{G\text{-}\mr{FC}}_{(S, \omega_S)}
\end{align}
the set of  
 $G$-FC quantizations on $(S, \omega_S)$.
We obtain the natural forgetting map
\begin{align} \label{W400}
\mfQ^{G\text{-}\mr{FC}}_{(S, \omega_S)} \longmigi \mfQ^{\mr{FC}}_{(S, \omega_S)}.
\end{align}

Furthermore, let $(T, \omega_{T})$ be the quotient of $(S, \omega_S)$ by the $G$-action.
The quotient morphism $f : S \migi T$ is a Galois \'{e}tale covering with Galois group $G$ and
$f^*(\omega_{T}) = \omega_S$.
Hence, pulling-back via $\pi$ induces a map of sets
\begin{align} \label{W1004}
\mfQ_{(T, \omega_T)}^{\mr{FC}} \longmigi \mfQ_{(S, \omega_S)}^{\mr{FC}}.
\end{align}
If $\mcO^\hslash_{T}$ is an FC quantization on $(T, \omega_{T})$,
then the pull-back $f^*(\mcO^\hslash_T)$ has naturally a structure of   $G$-FC quantization
  since the $G$-actions on $S$ and  $S^{(1)}$
 are compatible via $F_{S/k}$.
Conversely, let $\mcO^\hslash_S$ be a $G$-FC quantization on $(S, \omega_S)$.
Then, the sheaf $f_*(\mcO^\hslash_S)^G$ of $G$-invariant sections of
 $f_*(\mcO^\hslash_S)$
 specifies an FC quantization on $(T, \omega_T)$.
  One verifies immediately that the assignments 
 $\mcO_T^\hslash \mapsto f^*(\mcO_T^\hslash)$ and $\mcO_S^\hslash \mapsto f_*(\mcO_S^\hslash)^G$
 give 
 a  bijection correspondence
\begin{align} \label{W1005}
\mfQ^{\mr{FC}}_{(T, \omega_T)} \stackrel{\sim}{\longmigi}\mfQ^{G\text{-}\mr{FC}}_{(S, \omega_S)} 
\end{align}
 making the following diagram commute:
\begin{align}
\vcenter{\xymatrix{
\mfQ^{\mr{FC}}_{(T, \omega_T)} \ar[rr]_{\sim}^{(\ref{W1005})} \ar[rd]_{(\ref{W1004})} & & \mfQ^{G\text{-}\mr{FC}}_{(S, \omega_S)} \ar[ld]^{(\ref{W400})}
\\
 & \mfQ^{\mr{FC}}_{(T, \omega_T)}.  &
 }}
\end{align}

\vspace{5mm}
\subsection{Formal Weyl algebras} \label{S15}
\leavevmode\\ \vspace{-4mm}

In this subsection, we recall a canonical (Frobenius-constant) quantization on  the affine space $\mbA^{2n} := \mr{Spec} (k[x_1, \cdots, x_n, y_1, \cdots, y_n])$ ($n>0$) 
 equipped with the symplectic structure 
\begin{align}
\omega^\mr{Weyl} := \sum_{i=1}^n d y_i \wedge d x_i.
\end{align}
Here, notice that the Poisson bracket $\{ -, -\}_{\omega}$ associated with $\omega^\mr{Weyl}$  is given by 
  $\{ f, g\}_{\omega} = \sum_{i=1}^n \frac{\partial f}{\partial y_i} \cdot \frac{\partial g}{\partial x_i} - \frac{\partial f}{\partial x_i} \cdot \frac{\partial g}{\partial y_i}$ (for any local sections $f, g \in \mcO_{\mbA^{2n}}$).

For each commutative ring $R$ over $k$, we define $W^{2n}_R$ to be the (noncommutative) $R[[\hslash ]]$-algebra $W^{2n} := R[x_1, \cdots, x_n, y_1, \cdots, y_n][[\hslash ]]$ equipped with the multiplication ``$*$" given by 
\begin{align}
f * g := \sum_{\alpha \in \mbZ_{\geq 0}^n} \frac{\hslash^{|\alpha|}}{\alpha !} \cdot \partial^\alpha_{\vec{y}}(f) \cdot  \partial^\alpha_{\vec{x}}(g)
\end{align}
 for any $f, g \in R[x_1, \cdots, x_n, y_1, \cdots, y_n]$.
Hence, the $R[[\hslash]]$-algebra $W^{2n}_R$ is generated by the elements $x_1, \cdots, x_n, y_1, \cdots, y_n$ subject to  relations
\begin{align}
[x_i, x_j] = [y_i, y_j] =0, \hspace{5mm} [y_i, x_j] = \delta_{ij} \cdot \hslash
\end{align}
for all $0 < i, j \leq n$.
Since the center of this algebra coincides with $R[x^p_1, \cdots x_n^p, y^p_1, \cdots, y_n^p][[\hslash ]]$, $W^{2n}_R$ may be thought of as an $R[x^p_1, \cdots x_n^p, y^p_1, \cdots, y_n^p][[\hslash ]]$-algebra.

Here, write $\mr{Sp}_{2n}$ for the symplectic group over $k$ of rank $n$, i.e.,   
\begin{align}
\mr{Sp}_{2n}(R) = \left\{ A \in \mr{GL}_{2n}(R) \, | \, {^t A} J_{2n} A = J_{2n}\right\}
\end{align}
 for each commutative ring $R$ over $k$, where $J_{2n}:= \begin{pmatrix} O & E \\ -E & O\end{pmatrix}$  ($E$ denotes the unit matrix of size $n$).
Each $A \in \mr{Sp}_{2n}(R)$ yields an automorphism $\eta_A$ of $W^{2n}_R$ given by $(x_1, \cdots, x_n, y_1, \cdots, y_n) \mapsto (x_1, \cdots, x_n, y_1, \cdots, y_n){^t A}$.
If $\mr{Aut}(W^{2n}_R)$ denotes the automorphism group of the $R[[\hslash ]]$-algebra $W^{2n}_R$, then  the assignment $A \mapsto \eta_A$ determines  an injective homomorphism
\begin{align} \label{E4522}
\mr{Sp}_{2n}(R) \migi \mr{Aut}(W^{2n}_R).
\end{align}

Notice that $W^{2n}_k$ gives rise to 
an $\mcO_{(\mbA^{2n})^{(1)}}[[\hslash]]$-algebra
\begin{align} \label{Efgklo}
\mcW^{2n}_k,
\end{align}
which specifies 
  an FC quantization on $(\mbA^{2n}, \omega^\mr{Weyl})$, as well as on $(\mbA^{2n\times}, \omega^\mr{Weyl})$ via restriction.
The variety $\mbA^{2n\times}$ admits 
a free action of  $\mu_2 := \{ \pm1\}$ such that the automorphism corresponding to $-1 \in \mu_2$ is given by $(x, y) \mapsto (-x, -y)$.
This action preserves $\omega^\mr{Weyl}$, and we obtain the quotient symplectic variety
\begin{align}
(\mbA^{2n \times}_{/\mu_2}, \omega_{/\mu_2}^\mr{Weyl}).
\end{align}
One verifies that $\mcW^{2n}_k$ is a $\mu_2$-FC quantization, which descends to a FC quantization
\begin{align} \label{E4889}
\mcW^{2n}_{/\mu_2}
\end{align}
on $(\mbA^{2n \times}_{/\mu_2}, \omega_{/\mu_2}^\mr{Weyl})$.
(In our discussion, we will  use this quantization only in the case $n=1$.)

\vspace{10mm}
\section{Frobenius-projective structures and related objects} \vspace{3mm}

In this section, we review a positive characteristic analogue of a complex projective structure, called a Frobenius-projective structure.
Also, we
 discuss various bijective correspondences between Frobenius-projective structures on a curve
and some equivalent objects, i.e., dormant $\mr{PGL}_2$-opers, dormant $(\mr{SL}_2, \mbL)$-opers, and  projective connections with a full set of solutions.

\vspace{5mm}
\subsection{Frobenius-projective structures} \label{S1}
\leavevmode\\ \vspace{-4mm}

Let $n$ be a positive integer and denote by
 $\mr{PGL}_{n+1}$ the projective linear group over $k$ of rank $n+1$, which 
 is naturally identified with the automorphism group  of $\mbP^n$.
Given each algebraic group $G$ over $k$ and a smooth variety $S$, we denote by 
$G_{S}$ the sheaf of groups on $S$ represented by $G$.
Write
\begin{align}
G_{S}^\mr{F} := F^{-1}_{S/k} (G_{S^{(1)}}) \ \left(\subseteq G_{S} \right).
\end{align}
Also, we shall write
\begin{align}
\mcP^{\text{\'{e}t}}
\end{align}
for the sheaf of sets on $S$ that assigns, to each open subscheme  $U$ of $S$, the set of \'{e}tale morphisms $\phi : U \migi \mbP^n$ from $U$ to $\mbP^n$.
Each local section $\phi$ of $\mcP^{\text{\'{e}t}}$ may be regarded, by taking its graph, as a local section $\Gamma_\phi : U \migi U \times \mbP^n$ of the trivial $\mbP^n$-bundle $U \times \mbP^n \xrightarrow{\mr{pr}_1} U$. 
The sheaf $\mcP^{\text{\'{e}t}}$ has a $(\mr{PGL}_{n})_S^\mr{F}$-action described as follows.
Let $U$ be an open subscheme, $\phi : U \migi \mbP^n$ an element of  $\mcP^{\text{\'{e}t}}(U)$, and
$A$ an element of $(\mr{PGL}_{n})_S^\mr{F}(U) \ \left(\subseteq \mr{PGL}_n (U) \right)$.
Then,  one verifies immediately that the composite 
\begin{align}
A (\phi) : U \xrightarrow{\Gamma_\phi} U \times \mbP^n \xrightarrow{A} U \times \mbP^n \xrightarrow{\mr{pr}_2} \mbP^n
\end{align}
specifies  an element of $\mcP^{\text{\'{e}t}}(U)$.
The assignment $(A, \phi) \mapsto A  (\phi)$ defines a $(\mr{PGL}_{n})_S^\mr{F}$-action on $\mcP^{\text{\'{e}t}}$.

\vspace{3mm}
\bde[cf. ~\cite{Hos}, \S\,2, Definition 2.1 for the case $n=1$] \label{D01} \leavevmode\\
 \ \ \ We shall say that a subsheaf  $\mcS^\heartsuit$ of $\mcP^{\text{\'{e}t}}$ is a {\bf Frobenius-projective structure (of level $1$)}  on $S$ if $\mcS^\heartsuit$ is closed under the $(\mr{PGL}_{n})_S^\mr{F}$-action  on $\mcP^{\text{\'{e}t}}$, and  moreover, forms a $(\mr{PGL}_{n})_S^\mr{F}$-torsor on $S$ with respect to the resulting $(\mr{PGL}_{n})_S^\mr{F}$-action  on $\mcS^\heartsuit$.
     \ede
\vspace{3mm}

 We shall write
 \begin{align} \label{Efllo2}
  \mfP \mfS^\mr{F}_S
 \end{align}
 for the set of Frobenius-projective structures on $S$.

\vspace{5mm}
\subsection{Dormant indigenous bundles} \label{S1}
\leavevmode\\ \vspace{-4mm}

In what follows, let us fix a smooth curve $X$.
Recall from, e.g., ~\cite{Wak}, \S\,2, Definition 2.1 (i), that an {\bf indigenous bundle} (or, a $\mr{PGL}_2$-oper) on $X$  is 
a triple $\mcE^\spadesuit :=(\mcE, \nabla_\mcE, \sigma_\mcE)$ consisting of 
a flat $\mbP^1$-bundle $(\mcE, \nabla_\mcE)$ on $X$
 (i.e., a pair of a $\mbP^1$-bundle  $\mcE$ on  $X$ and a connection $\nabla_\mcE$ on $\mcE$) 
  and  a global section $\sigma_\mcE : X \migi \mcE$ such that the Kodaira-Spencer map $\mcT_{X/k} \migi \sigma_\mcE^*(\mcT_{\mcE/X})$ associated to $\sigma_\mcE$ 
    is   nowhere vanishing.
   (We omit to describe in detail  the definition of an indigenous bundle because it will not be necessary for our discussion.)
We shall say that an indigenous bundle $\mcE^\spadesuit := (\mcE, \nabla_\mcE, \sigma_\mcE)$ is {\bf dormant} if the connection $\nabla_\mcE$ has vanishing $p$-curvature (cf. ~\cite{Wak}, \S\,3, Definition 3.1).
Write
\begin{align} \label{E69207}
\mfI \mfB_X \left( \text{resp.,} \ \mfI \mfB_X^{^\mr{Zzz...}}\right)
\end{align}
for the set of isomorphism classes of indigenous bundles (resp., dormant indigenous bundles) on $X$.
Then, 
there exists a canonical bijection of sets
\begin{align} \label{W101}
\mfP \mfS^\mr{F}_X \isom \mfI \mfB_X^{^\mr{Zzz...}}.
\end{align}
(~\cite{Hos}, \S\,3, Proposition 3.11 in the case where $X$ is proper).
In fact, let $\mcS^\heartsuit$ be a Frobenius-projective structure on $X$.
The $\mr{PGL}_2$-torsor over $X^{(1)}$ corresponding to $\mcS^\heartsuit$ via
the underlying homeomorphism of $F_{X/k}$ specifies  
 a $\mbP^1$-bundle 
  over $X^{(1)}$.
The pull-back $\mcE_\mcS$ of this $\mbP^1$-bundle over $X$ 
 admits naturally a connection $\nabla^\mr{can}_{\mcE_\mcS}$ with vanishing $p$-curvature (cf. ~\cite{Kat}, \S\,5, Theorem 5.1).
Moreover, the local sections $U \migi U \times \mbP^1$ (for various open subschemes $U$ of $X$) classified by sections in $\mcS^\heartsuit$ may be glued together to obtain a well-defined global section $\sigma_{\mcE_\mcS} : X \migi \mcE_\mcS$. It follows from the definition of a Frobenius-projective structure that
the resulting triple 
\begin{align} \label{W60}
\mcE_\mcS^\spadesuit := (\mcE_\mcS, \nabla^\mr{can}_{\mcE_\mcS}, \sigma_{\mcE_\mcS})
\end{align}
 specifies a dormant indigenous bundle on $X$.
The resulting  assignment $\mcS^\heartsuit \mapsto \mcE_\mcS^\spadesuit$ gives  the bijection (\ref{W101}).

\vspace{3mm}
\begin{rema} \label{R111}
\leavevmode\\
\ \ \
If $X$ is proper, then we know an explicit formula for computing the number of dormant indigenous bundles on  $X$, as proved in  
~\cite{Wak}, Theorem A. 
In particular, there exists at least one dormant indigenous bundle on any (not necessarily proper) smooth curve as it has a smooth compactification.
\end{rema}

\vspace{5mm}
\subsection{Dormant indigenous vector  bundles} \label{S1}
\leavevmode\\ \vspace{-4mm}

In this subsection, we describe indigenous bundles  and their higher-rank generalizations in terms of vector  bundles.
Let $n$ be a positive integer with $n<p$.
Here, recall that, for each vector bundle $\mcF$ on $X$ of rank $n$,
a {\it connection} on  $\mcF$ means a $k$-linear morphism $\nabla_\mcF : \mcF \migi \Omega_{X}\otimes \mcF$ satisfying that $\nabla_\mcF (a \cdot v) = da \otimes v + a \cdot \nabla_\mcF (v)$ for any local sections $a \in \mcO_X$ and $v \in \mcF$. 
Given such a connection $\nabla_\mcF$, we have a connection $\nabla_{\mr{det}(\mcF)}$ on the determinant bundle $\mr{det}(\mcF)$ induced by $\nabla_\mcF$, i.e., given by $\nabla_{\mr{det}(\mcF)}(a_{1} \wedge \cdots \wedge a_{n}) = \sum_{i=1}^n a_1 \wedge \cdots \wedge \nabla_\mcF (a_i) \wedge \cdots \wedge a_n$,  where $n:=\mr{rk}(\mcF)$.

Recall (cf. ~\cite{BD2}, \S\,2.1) that a {\bf $\mr{GL}_n$-oper} on $X$ is a collection of data
\begin{align}
(\mcF, \nabla_\mcF, \mcF^\bullet)
\end{align}
consisting of a rank $n$ vector bundle $\mcF$ on $X$, a connection $\nabla_\mcF$
  on $\mcF$, and an $n$-step decreasing filtration $\mcF^\bullet := \{ \mcF^j \}_{j=0}^n$ on $\mcF$ satisfying the following conditions
\begin{itemize}
\item[$\bullet$]
Each $\mcF^j$  is a subbundle of $\mcF$ such that
$\mcF^0 = \mcF$, $\mcF^n = \mcF$, and $\mr{gr}^j_\mcF := \mcF^j/\mcF^{j+1}$ ($0\leq j \leq n-1$) is a line bundle;
\item[$\bullet$]
$\nabla_\mcF (\mcF^j) \subseteq \Omega_X \otimes \mcF^{j-1}$ ($1 \leq j \leq n-1$) and the morphism
$\mr{KS}^j_{\mcF^\bullet}  : 
\mr{gr}^j_\mcF \migi \Omega_X \otimes \mr{gr}^{j-1}_\mcF$
induced by $\nabla_\mcF$ (which is verified to be $\mcO_X$-linear) is an isomorphism.
\end{itemize}
In particular, 
 a $\mr{GL}_2$-oper
  on $X$ is determined by  a triple $(\mcF, \nabla_\mcF, \mcN_\mcF)$ consisting of a pair $(\mcF, \nabla_\mcF)$  as above (with $n=2$) and a line subbundle  $\mcN_\mcF$ of $\mcF$ such that
 the $\mcO_X$-linear composite
  \begin{align} \label{W80}
 \mr{KS}_{\mcN_\mcF} : \mcN_\mcF \migiincl \mcF \xrightarrow{\nabla_\mcF} \Omega_{X}\otimes \mcF \migisurj \Omega_X \otimes (\mcF/\mcN_\mcF)
 \end{align}
(called the {\it Kodaira-Spencer map} associated with $\mcN_\mcF$) is an isomorphism.

Next, let  us fix a theta characteristic $\mbL := (\mcL, \psi_\mcL)$ (cf. \S\,\ref{S8}) of $X$.

\vspace{3mm}
\bde[cf. ~\cite{Wak}, \S\,2, Definition 2.3 (ii), in the case $n=2$] \label{D10128} \leavevmode\\
 \ \ \ 
 An {\bf $(\mr{SL}_n, \mbL)$-oper} on $X$
is a collection of data
\begin{align}
\mcF^\diamondsuit :=(\mcF, \nabla_\mcF, \mcF^\bullet, \eta_\mcF),
\end{align}
where $(\mcF, \nabla_\mcF, \mcF^\bullet)$ is a $\mr{GL}_n$-oper on $X$ and $\eta_\mcF$ denotes an isomorphism $\mcL^{\otimes (n-1)} \isom \mcF^{n-1}$
such that 
  the connection $\nabla_{\mr{det}(\mcF)}$ on $\mr{det}(\mcF)$ 
   coincides with $d : \mcO_X \migi \Omega_X$ via the composite  isomorphism 
 \begin{align} \label{Effgh3}
\delta_{\mcF^\diamondsuit} :  \mr{det}(\mcF)\left(= \bigotimes_{j=0}^{n-1} \mr{gr}^j_\mcF \right) &\xrightarrow{\otimes_j \kappa_j} \bigotimes_{j=0}^{n-1} \Omega_X^{\otimes (j-n+1)} \otimes \mcF^{n-1} \\
 & \xrightarrow{\sim} \Omega_{X}^{\otimes \frac{-n(n-1)}{2}} \otimes (\mcF^{n-1})^{\otimes n} \notag \\
 & \xrightarrow{\sim} \Omega_{X}^{\otimes \frac{-n(n-1)}{2}} \otimes  \mcL^{\otimes n (n-1)}\notag \\
 & \xrightarrow{\sim} \mcO_X, \notag
  \end{align}
 where the third isomorphism arises from $\eta_\mcF$, the fourth isomorphism arises from $\psi_\mcF$, and 
each $\kappa_j$ ($j=0, \cdots, n-1$) in the first isomorphism denotes the composite isomorphism
\begin{align}
\mr{gr}^j_\mcF \isom \Omega^{\otimes (-1)}_X \otimes \mr{gr}^{j+1}_\mcF \isom \cdots \isom \Omega_{X}^{\otimes (j-n+1)} \otimes \mr{gr}^{n-1}_\mcF \ \left(=\Omega_{X}^{\otimes (j-n+1)} \otimes \mcF^{n-1} \right)
\end{align}
 each of whose constituent arises from $\mr{KS}^{(-)}_{\mcF^\bullet}$.
In a natural manner, we can define the notion of an isomorphism between $(\mr{SL}_n, \mbL)$-opers.
 Also, we shall say that an $(\mr{SL}_n, \mbL)$-oper is {\bf dormant} if it has vanishing $p$-curvature.
\ede
\vspace{3mm}

We denote by
\begin{align}
\mfO \mfp_{(\mr{SL}_n, \mbL), X} \left(\text{resp.,} \ \mfO \mfp_{(\mr{SL}_n, \mbL), X}^{^\mr{Zzz...}}\right)
\end{align}
the set of isomorphism classes of $(\mr{SL}_n, \mbL)$-opers (resp., dormant $(\mr{SL}_n, \mbL)$-opers) on $X$.
According to  ~\cite{Wak}, \S\,2, Proposition 2.4,  there exists a canonical bijection
\begin{align} \label{E02f34}
\mfI \mfB_{X} \isom \mfO \mfp_{(\mr{SL}_n, \mbL), X},
\end{align}
which restricts to a bijection
\begin{align} \label{E02234}
\mfI \mfB_{X}^{^\mr{Zzz...}} \isom \mfO \mfp_{(\mr{SL}_n, \mbL), X}^{^\mr{Zzz...}}.
\end{align}
Let  $\mcE^\spadesuit := (\mcE, \nabla_\mcE, \sigma_\mcE)$ be an indigenous bundle on $X$ and $\mcF^\diamondsuit := (\mcF, \nabla_\mcF, \mcN_\mcF, \eta_\mcF)$ denote the $(\mr{SL}_2, \mbL)$-oper corresponding to $\mcE^\spadesuit$.
Then, $(\mcE, \nabla_\mcE)$ may be obtained from $(\mcF, \nabla_\mcF)$ via  projectivization, and $\eta_\mcF : \mcL \isom \mcN_\mcF \ \left(\subseteq \mcF \right)$ induces a morphism
\begin{align} \label{Ega3t}
\eta^\A_\mcF : \A (\mcL)^\times \migi \A (\mcF)^\times
\end{align}
over $X$ 
such that the following square diagram 
 \begin{align} \label{Ddd01}
\vcenter{\xymatrix@C=36pt@R=36pt{
 \A (\mcL)^\times \ar[r]^{\eta^\A_\mcF} \ar[d]  
  & \A (\mcF)^\times   \ar[d]^{\pi_\mcF}
 \\
 X \ar[r]_-{\sigma_\mcE}  &  \mcE \ \left(= \mbP (\mcF) \right)
 }}
 \end{align}
is commutative and cartesian.

\vspace{3mm}
\begin{exa} \label{E510d}
\leavevmode\\
\ \ \
Let $\mcF^\diamondsuit := (\mcF, \nabla_\mcF, \mcN_\mcF, \eta_\mcF)$ be an
$(\mr{SL}_2, \mbL)$-oper
 on $X$.
The $(n-1)$-th symmetric power $\mcF_{\mr{SL}_n} := S^{n-1}(\mcF)$ of  $\mcF$ is a rank $n$ vector bundle.
For each $j=0, 1, \cdots, n$, the image  $\mcF^j_{\mr{SL}_n}$
of the natural morphism $\mcN_\mcF^{\otimes j} \otimes S^{n-1-j}(\mcF)\migi S^{n-1}(\mcF)$ is 
  a rank $(n-j)$ subbundle of $\mcF_{\mr{SL}_n}$.
The collection $\mcF_{\mr{SL}_n}^\bullet := \{ \mcF_{\mr{SL}_n}^j \}_{j=0}^n$ forms an $n$-step decreasing filtration on $\mcF_{\mr{SL}_n}$.
Let $\eta_{\mcF, \mr{SL}_n} : \mcL^{\otimes (n-1)} \isom \mcF_{\mr{SL}_n}^{n-1}$ be the composite isomorphism of $\eta_\mcF^{\otimes (n-1)} : \mcL^{\otimes (n-1)} \isom \mcN_\mcF^{\otimes (n-1)}$ with the natural isomorphism $\mcN^{\otimes (n-1)}_{\mcF} \isom \mcF_{\mr{SL}_n}^{n-1}$.
Also, let $\nabla_{\mcF, \mr{SL}_n}$ be the  connection  on $\mcF_{\mr{SL}_n}$ induced naturally by $\nabla_\mcF$.
Then, one verifies immediately that the collection of data
\begin{align}
\mcF^\diamondsuit_{\mr{SL}_n} := (\mcF_{\mr{SL}_n}, \nabla_{\mcF, \mr{SL}_n}, \mcF^\bullet_{\mr{SL}_n}, \eta_{\mcF, \mr{SL}_n})
\end{align}
forms an $(\mr{SL}_n, \mbL)$-oper on $X$.
If, moreover, $\mcF^\diamondsuit$ is dormant, then the resulting $(\mr{SL}_n, \mbL)$-oper  $\mcF^\diamondsuit_{\mr{SL}_n}$ is dormant.
Thus, the assignment $\mcF^\diamondsuit \mapsto \mcF^\diamondsuit_{\mr{SL}_n}$ defines a map of sets
\begin{align} \label{Erty7}
\mfO \mfp_{(\mr{SL}_2, \mbL), X} \migi  \mfO \mfp_{(\mr{SL}_n, \mbL), X}
\end{align}
which restricts to a map  $\mfO \mfp_{(\mr{SL}_2, \mbL), X}^{^\mr{Zzz...}} \migi  \mfO \mfp_{(\mr{SL}_n, \mbL), X}^{^\mr{Zzz...}}$.
\end{exa}

\vspace{5mm}
\subsection{Projective connections} \label{S63}
\leavevmode\\ \vspace{-4mm}

In this subsection, we discuss higher-order  projective connections.
Let $\mbL$ be as above and let $D : \mcL^{\otimes (-n+1)} \migi \mcL^{\otimes (n+1)}$ be an $n$-th order  differential operator 
(i.e., an element of $\mcD {\it iff}^{\leq n}_{ \mcL^{\otimes (-n+1)}, \mcL^{\otimes (n+1)}}$) satisfying the equality $\Sigma (D)=1$ (cf. \S\,\ref{S83}) under the identification
\begin{align}
\mcH om_{\mcO_X}(\mcL^{\otimes (-n+1)}, \mcL^{\otimes (n+1)}\otimes S^n(\mcT_X)) \isom  \mcH om_{\mcO_X} (\mcL^{\otimes (-n+1)}, \mcL^{\otimes (-n+1)})  \isom \mcO_X
\end{align}
induced by $\psi_\mcL$.
Denote by ${^t D}$ the transpose of $D$, which is a differential operator
$\Omega_X \otimes (\mcL^{\otimes (n+1)})^\vee \migi \Omega_X \otimes (\mcL^{\otimes (-n+1)})^\vee)$.
 If $D$ is locally expressed as $D = \sum_{i=0}^n a_i \cdot \partial^i$ (for a local generator $\partial$ of $\mcT_X$), then
${^t D} = \sum_{i=0}^n (-\partial)^i \cdot a_i$.
Since $\psi_\mcL$ allows us to consider  $\Omega_X \otimes (\mcL^{\otimes (n+1)})^\vee$ and  $\Omega_X \otimes (\mcL^{\otimes (-n+1)})^\vee)$  as $\mcL^{\otimes (-n+1)}$ and $\mcL^{\otimes (n+1)}$  respectively, $D$ may be thought of as a differential operator in   $\mcD {\it iff}^{\leq n}_{ \mcL^{\otimes (-n+1)}, \mcL^{\otimes (n+1)}}$.
In particular, it makes sense to speak of   the operator  $D' := \frac{1}{2} \cdot (D -(-1)^{n}  \cdot {^t D})$.
Moreover, by the equality $\Sigma (D) =1$, the operator $D'$ turns out to be of order $n-1$.
We refer to the principal symbol $\Sigma_\mr{sub}(D) := \Sigma (D')$ of $D'$ as the {\it subprincipal symbol} of $D$.

\vspace{3mm}
\bde \label{D01ff} \leavevmode\\
 \ \ \ 
 An {\bf $n$-th order  projective connection} for $\mbL$ is 
 an $n$-th order differential operator $D^\clubsuit : \mcL^{\otimes (-n+1)} \migi \mcL^{\otimes (n+1)}$ with $\Sigma (D^\clubsuit) =1$ and $\Sigma_\mr{sub}(D^\clubsuit) =0$.
For simplicity, we refer to each second order projective connection as a {\bf projective connection}.
      \ede
\vspace{3mm}

We shall write
\begin{align} \label{E41180}
\mfP \mfC^n_{X, \mbL} \ \left(\text{resp.,} \ \mfP \mfC_{X, \mbL}^{n, \mr{full}} \right)
\end{align}
for 
the set of $n$-th order  projective connections for $\mbL$ (resp., the set of $n$-th order projective connections for $\mbL$ having a full set of solutions).

\vspace{3mm}
\bpr[cf. ~\cite{BD2}, \S\,2.1 and \S\,2.8; ~\cite{Kat2}, Proposition 6.0.5] \label{L100} \leavevmode\\
 \ \ \ 
There exists a canonical bijection 
\begin{align} \label{F001}
\mfO \mfp_{(\mr{SL}_n, \mbL), X} \isom \mfP \mfC^n_{X, \mbL}
\end{align}
 restricting to a bijection
 \begin{align} \label{F002}
 \mfO \mfp_{(\mr{SL}_n, \mbL), X}^{^\mr{Zzz...}} \isom \mfP \mfC^{n, \mr{full}}_{X, \mbL}.
 \end{align}
 \epr
\begin{proof}
First, we shall construct the bijection (\ref{F001}).
Let  $\mcF^{\diamondsuit}:= (\mcF, \nabla_\mcF, \mcF^\bullet, \eta_\mcF)$ be an $(\mr{SL}_n, \mbL)$-oper  on $X$.
The connection $\nabla_\mcF$ 
induces, inductively on  $i \leq n$, 
an $\mcO_X$-linear morphism
$\nabla_\mcF^{\mcD, i} : \mcD_X^{\leq i} \otimes \mcF \migi \mcF$ 
determined by the condition that $\nabla^{\mcD, 0}_\mcF = \mr{id}_\mcF$ and $\nabla_\mcF^{\mcD, i}(\partial^i \otimes v) = \langle \partial, \nabla_\mcF(\nabla_\mcF^{\mcD, i-1} (\partial^{i-1} \otimes v))\rangle$ for any local generator  $\partial \in \mcT_X$ and any local section $v \in \mcF$.
By   the definition of a $\mr{GL}_n$-oper, we see (by induction on $i$) that
the morphism $\nabla_\mcF^{\mcD, i}$ for $i \leq n-1$ restricts to an isomorphism $\mcD_{X}^{\leq i} \otimes \mcF^{n-1} \migi \mcF^{n-i-1}$ and hence $\nabla_\mcF^{\mcD, n}$ is surjective.
The composite $(\nabla^{\mcD, n-1}_{\mcF})^{-1} \circ \nabla^{\mcD, n}_{\mcF}$, regarded as an $\mcO_X$-linear morphism   $\mcD_X^{\leq n} \otimes \mcL^{\otimes (n-1)} \migi \mcD_X^{\leq (n-1)} \otimes \mcL^{\otimes (n-1)}$ via $\eta_\mcF$, 
determines a split surjection of the following short exact sequence:
 \begin{align} \label{Eddfg}
 0 \longmigi \mcD^{\leq (n-1)}_X \otimes \mcL^{\otimes (n-1)} \longmigi \mcD^{\leq n}_X \otimes \mcL^{\otimes (n-1)} \longmigi (\mcD^{\leq n}_X/\mcD^{\leq (n-1)}_X) \otimes \mcL^{\otimes (n-1)} \longmigi 0.
 \end{align}
Let us consider the composite $\mcL^{\otimes (-n-1)} \migi  \mcD_X^{\leq n} \otimes \mcL^{\otimes (n-1)}$ of the corresponding  split injection 
$\left(\mcT_X^{\otimes n} \otimes \mcL^{\otimes (n-1)}= \right) \ (\mcD^{\leq n}_X/\mcD^{\leq (n-1)}_X) \otimes \mcL^{\otimes (n-1)} \migiincl \mcD_X^{\leq n} \otimes \mcL^{\otimes (n-1)}$ 
and the isomorphism $\mcL^{\otimes (-n-1)} \isom \mcT_X^{\otimes n} \otimes \mcL^{\otimes (n-1)}$ induced naturally by $\psi_\mcL$;
it corresponds to an $\mcO_X$-linear morphism $\mcL^{\otimes (-n+1)} \migi \mcL^{\otimes (n+1)} \otimes \mcD_X^{\leq n}$, or equivalently, an $n$-th differential operator $D^\clubsuit_{\mcF} : \mcL^{\otimes (-n+1)} \migi \mcL^{\otimes (n+1)}$ (cf. (\ref{Edr56})). One verifies immediately that $\Sigma (D^\clubsuit_{\mcF}) =1$, and moreover, (by taking account of the fact that $(\mr{det}(\mcF), \nabla_{\mr{det}(\mcF)}) \cong (\mcO, d)$) that $\Sigma_\mr{sub}(D^\clubsuit_{\mcF}) =0$.
Thus, $D^\clubsuit_{\mcF}$ specifies a projective connection for $\mbL$.
Thus, the assignment $\mcF^\diamondsuit \mapsto D^\clubsuit_{\mcF}$  defines a map of sets $\mfO \mfp_{(\mr{SL}_n, \mbL), X} \migi \mfP \mfC^n_{X, \mbL}$.

Conversely, let $D^\clubsuit$ be a projective connection belonging to  $\mfP \mfC^n_{X, \mbL}$.
For each $i=0, \cdots, n$, we shall write
$\mcF_D^i := \mcD_{X}^{\leq (n-i-1)}\otimes \mcL^{\otimes (n-1)}$ and $\mcF_D^\bullet := \{ \mcF_D^i \}_{i=0}^n$.
The operator $D^\clubsuit$ may be thought of as an $\mcO_X$-linear morphism $\mcL^{\otimes (-n+1)} \migi \mcL^{\otimes (n+1)} \otimes \mcD_{X}^{\leq n}$ (via (\ref{Edr56})), or equivalently, $\left((\mcD_X^{\leq n}/\mcD_X^{\leq (n-1)}) \otimes \mcL^{\otimes (n-1)}= \right) \mcL^{\otimes (-n-1)} \migi \mcD_X^{\leq n} \otimes \mcL^{\otimes (n-1)}$.
It specifies a split injection of (\ref{Eddfg}),
where we shall write $\nabla'$ for the corresponding split surjection $\mcD_X^{\leq n} \otimes \mcL^{\otimes (n-1)} \migisurj \mcD_X^{\leq (n-1)} \otimes \mcL^{\otimes (n-1)} \ \left(= \mcF^0_D \right)$.
Then, there exists a unique connection $\nabla_D$ on $\mcF_D^0 \ \left(\subseteq \mcD_X^{\leq n}\otimes \mcL^{\otimes (n-1)} \right)$ determined by the condition  that
$\langle \partial, \nabla_D (\partial^i \otimes v)\rangle = \nabla' (\partial^{i+1} \otimes v)$ ($i= 0, \cdots, n-1$) for any local generator  $\partial \in \mcT_X$ and any local section $v \in \mcL^{\otimes (n-1)}$.
  If $\eta_D$ denotes the natural isomorphism $\mcL^{\otimes (n-1)} \isom \mcF_D^{n-1}$, then (because of the assumption that $\Sigma (D) =1$ and $\Sigma_\mr{sub}(D)=0$) the collection 
  $\mcF_D^\diamondsuit := (\mcF^0_D, \nabla_D, \mcF_D^\bullet, \eta_D)$ forms an  $(\mr{SL}_n, \mbL)$-oper on $X$.
  One verifies that the assignment $D^\clubsuit \mapsto \mcF_D^\diamondsuit$ turns out to be the inverse of the map $\mfO \mfp_{(\mr{SL}_n, \mbL), X} \migi \mfP \mfC^n_{X, \mbL}$ obtained above, which completes the former assertion.
 
 Next, we shall  consider the latter assertion.
 Let us take  a projective connection $D^\clubsuit$ belonging to  $\mfP \mfC^n_{X, \mbL}$, and denote by $\mcF^\diamondsuit := (\mcF, \nabla_\mcF, \mcF^\bullet, \eta_\mcF)$ the corresponding $(\mr{SL}_n, \mbL)$-oper  constructed by the above steps.
 If $D^\clubsuit$ may be expressed (after choosing a local identification $\mcL \cong \mcO_X$) locally  as $D^\clubsuit = \partial^n + q_1 \partial^{n-1}+ \cdots + q_{n-1}\partial + q_n$ (for a local generator $\partial \in \mcT_X$ and local sections $q_1, \cdots, q_n \in \mcO_X$),
then the dual connection $\nabla_\mcF^\vee$ of $\nabla_\mcF$
may be expressed locally (with respect to a suitable local basis) as
\begin{align}
\nabla_D^\vee = \partial - \begin{pmatrix} -q_1 & -q_2& -q_3& \cdots & -q_{n-1}& - q_n \\ 1& 0 & 0& \cdots & 0 & \\ 0 & 1 & 0& \cdots & 0 & 0 \\ 0 & 0 &1& \cdots & 0& 0 \\ \vdots & \vdots & \vdots & \ddots& \vdots & \vdots  \\ 0 & 0 & 0 & \cdots & 1 & 0  \end{pmatrix}.
\end{align}
Then, $y \mapsto {^t(}\partial^{n-1}(y), \cdots, \partial (y), y)$
gives a bijective correspondence between the solutions of the differential equation $D^\clubsuit(y) =0$ and the horizontal (with respect to $\nabla_\mcF^\vee$) local sections of $\mcF^\vee$.
This implies that $D^\clubsuit$ has a full set of solutions if and only if the connection $\nabla_\mcF^\vee$, as well as $\nabla_\mcF$, has vanishing $p$-curvature.
Consequently, the  bijection $\mfO \mfp_{(\mr{SL}_n, \mbL), X} \isom \mfP \mfC^n_{X, \mbL}$ restricts to a bijection $\mfO \mfp^{^\mr{Zzz...}}_{(\mr{SL}_n, \mbL), X} \isom \mfP \mfC^{n, \mr{full}}_{X, \mbL}$, as desired.
\end{proof}
\vspace{3mm}

Thus, we have obtained  various maps  of sets, as displayed below:
\begin{align} \label{E46090}
\vcenter{\xymatrix{
\mfP \mfS_X^\mr{F} \ar[r]_-\sim^-{(\ref{W101})} & \mfI \mfB^{^\mr{Zzz...}}_X \ar[r]_-\sim^-{(\ref{E02234})} \ar[d]^{\mr{incl.}} & \mfO \mfp^{^\mr{Zzz...}}_{(\mr{SL}_2, \mbL), X} \ar[r]^-{(\ref{F002})}_-\sim \ar[d]^{\mr{incl.}} & \mfP \mfC^{2, \mr{full}}_{X, \mbL}  \ar[d]^{\mr{incl.}}\\
& \mfI \mfB_X \ar[r]^-\sim_-{(\ref{E02f34})} & \mfO  \mfp_{(\mr{SL}_2, \mbL), X} \ar[r]^-\sim_-{(\ref{F001})} &  \mfP \mfC^2_{X, \mbL},
}}
\end{align}
where
 all the vertical arrows are natural inclusions.
Moreover, there is a map for the set of projective connections (resp., having a full set of solutions) to the set of $n$-th order projective connections (resp., having a full set of solutions)
\begin{align} \label{E4289}
\xi^{2 \migi n}  : \mfP \mfC^2_{X, \mbL} \migi \mfP \mfC^n_{X, \mbL} \ \left(\text{resp.,} \ \xi^{2 \migi n} : \mfP \mfC^{2, \mr{full}}_{X, \mbL} \migi \mfP \mfC^{n, \mr{full}}_{X, \mbL} \right)
\end{align}
constructed in such a way that the  square diagram
\begin{align}
\vcenter{\xymatrix{
\mfO \mfp_{(\mr{SL}_2, \mbL), X} \ar[r]^{(\ref{Erty7})} \ar[d]^\wr_{(\ref{F001})} & \mfO \mfp_{(\mr{SL}_n, \mbL), X} \ar[d]_\wr^{(\ref{F001})}
\\
\mfP \mfC^2_{X, \mbL} \ar[r]_{\xi^{2 \migi n}} & \mfP \mfC^n_{X, \mbL}
}}
\
\left(\text{resp.,} \vcenter{\xymatrix{\mfO \mfp_{(\mr{SL}_2, \mbL), X} \ar[r]^{(\ref{Erty7})} \ar[d]^\wr_{(\ref{F001})} & \mfO \mfp_{(\mr{SL}_n, \mbL), X} \ar[d]_\wr^{(\ref{F001})}
\\
\mfP \mfC^2_{X, \mbL} \ar[r]_{\xi^{2 \migi n}} & \mfP \mfC^n_{X, \mbL}
}} \right)
\end{align}
is commutative.

\vspace{5mm}
\subsection{Case of the projective line} \label{S1}
\leavevmode\\ \vspace{-4mm}

In this subsection, we shall consider the case where 
$X= \mbP^1 \ \left(= \mr{Proj}(k [x, y]) \right)$ equipped with the theta characteristic $\mbL_0$ (cf. (\ref{E09076})).
We will observe that,  in this case, there is a typical example of a Frobenius-projective structure (resp., a dormant indigenous bundle; resp., a dormant $(\mr{SL}_n, \mbL_0)$-oper; resp., a projective connection for $\mbL_0$ having a full set of solutions) on $\mbP^1$, which  will be denoted by $\mcS_0^\heartsuit$ (resp., $\mcE^\spadesuit_0$; resp., $\mcF^{\diamondsuit}_0$; resp., $D_0^\clubsuit$).

First, we define
\begin{align}
\mcS_0^\heartsuit
\end{align}
to be the subsheaf of $\mcP^{\text{\'{e}t}}$ (for $X = \mbP^1$) which, to any open subscheme $U$ of $\mbP^1$, assigns the set  
\begin{align}
S_0^\heartsuit(U) := 
\{ A (\mr{op}_U) \in \mcP^{\text{\'{e}t}}(U) \, | \, A \in (\mr{PGL}_{2})_{\mbP^1}^\mr{F} (U) \},
\end{align}
 where $\mr{op}_U$ denotes the natural open immersion $U \migiincl \mbP^1$.
Then, $\mcS_0^\heartsuit$ forms a trivial $(\mr{PGL}_2)_{\mbP^1}^\mr{F}$-torsor, and hence,  specifies a Frobenius-projective structure on $\mbP^1$.

Next, we shall write $\mcE_0 := \mbP^1 \times \mbP^1$, which defines 
the trivial $\mbP^1$-bundle on $\mbP^1$ by 
regarding 
  the first projection $\mr{pr}_1 : \mcE_0 \migi \mbP^1$ 
  as its structure morphism.
Write $\nabla_0$ for the trivial connection on this trivial $\mbP^1$-bundle; it is clear that $\nabla_0$  has vanishing $p$-curvature.
The Kodaira-Spencer map (with respect to $\nabla_0$) of the diagonal embedding $\sigma_0 : \mbP^1 \migi \mbP^1 \times \mbP^1 \ \left(= \mcE_0 \right)$ is nowhere vanishing.
Thus, the triple
\begin{align} \label{E469}
\mcE^\spadesuit_0 := (\mcE_0, \nabla_0, \sigma_0)
\end{align}
 forms a dormant indigenous bundle on $\mbP^1$.

Moreover, recall the injection $\eta_0: \mcO_{\mbP^1}(-1) \migiincl \mcO_{\mbP^1}^{\oplus 2}$  introduced in Example \ref{E5d10}, (\ref{Effgh}) (of  the case $n=1$), which we shall
identify with the resulting  isomorphism from $\mcO_{\mbP^1}(-1)$ onto its image.
Then, the collection
\begin{align}
\mcF^\diamondsuit_0 := (\mcO_{\mbP^1}^{\oplus 2}, d^{\oplus 2}, \mcO_{\mbP^1}(-1), \eta_{0})
\end{align}
forms a dormant $(\mr{SL}_2, \mbL_{0})$-oper on $\mbP^1$.

Finally, 
let us consider the second order differential operator $D_{0, u}^\clubsuit : \mcO_{\mbP^1}(1)|_U \migi \mcO_{\mbP^1}(-3)|_U$ (resp., $D_{0, t}^\clubsuit : \mcO_{\mbP^1}(1)|_T \migi \mcO_{\mbP^1}(-3)|_T$) on the open subscheme $U := \mr{Spec}(k[u])$ (resp., $T := \mr{Spec}(k[t])$) of $\mbP^1$, where $u := x/y$ (resp., $t:= y/x$), given by 
$f (u) \cdot y \mapsto  \frac{\partial^2 f}{\partial u^2}(u) \cdot y^{-3}$
(resp., $g (t) \cdot x \mapsto  \frac{\partial^2 g}{\partial t^2}(t) \cdot x^{-3}$).
Then, 
$D_{0, u}^\clubsuit$ and $D_{0, t}^\clubsuit$ may be glued together to obtain a globally defined differential operator
\begin{align}
D_0^\clubsuit : \mcO_{\mbP^1}(-1)^\vee \ \left(=\mcO_{\mbP^1}(1) \right) \migi \mcO_{\mbP^1}(-1)^{\otimes 3} \ \left(= \mcO_{\mbP^1}(-3) \right)
\end{align}
forming a  projective connection for $\mbL_0$.

\vspace{3mm}
\bpr \label{P205} \leavevmode\\
 \ \ \ 
All the sets  $\mfP \mfS_{\mbP^1}^\mr{F}$, $\mfI \mfB_{\mbP^1}$, $\mfI \mfB_{\mbP^1}^{^\mr{Zzz...}}$, $\mfO \mfp_{(\mr{SL}_2, \mbL_0), \mbP^1}$,  $\mfO \mfp_{(\mr{SL}_2, \mbL_0), \mbP^1}^{^\mr{Zzz...}}$, $\mfP \mfC_{\mbP^1, \mbL_0}^2$, and $\mfP \mfC_{\mbP^1, \mbL_0}^{2, \mr{full}}$ are singletons respectively.
That is to say,
\begin{align}
\mfP \mfS_{\mbP^1}^\mr{F} &= \left\{\mcS^\heartsuit_0 \right\}, \\
\mfI \mfB_{\mbP^1}^{^\mr{Zzz...}} & =\mfI \mfB_{\mbP^1} = \left\{ \mcE^{\spadesuit}_0\right\},  \notag\\
\mfO \mfp_{(\mr{SL}_2, \mbL_0), \mbP^1} & = \mfO \mfp_{(\mr{SL}_2, \mbL_0), \mbP^1}^{^\mr{Zzz...}} = \left\{ \mcF^{\diamondsuit}_0\right\}, \notag \\
\mfP \mfC^2_{\mbP^1, \mbL_0} & = \mfP \mfC^{2, \mr{full}}_{\mbP^1, \mbL_0} = \left\{ D^\clubsuit_0 \right\}. \notag
\end{align}

 \epr
\begin{proof}
By the various bijections in (\ref{E46090}), 
it suffices to prove that $\mfP \mfC^2_{\mbP^1, \mbL_0}$ contains at most one element.
As discussed in the proof of Proposition \ref{L100},
each element of $\mfP \mfC^2_{\mbP^1, \mbL_0}$ corresponds to
a splitting of 
the short exact sequence
\begin{align} \label{E009ju}
0 \longmigi \mcD_{\mbP^1}^{\leq 1} \otimes \mcO_{\mbP^1}(-1) \longmigi
\mcD_{\mbP^1}^{\leq 2} \otimes \mcO_{\mbP^1}(-1) \longmigi
(\mcD_{\mbP^1}^{\leq 2} /\mcD_{\mbP^1}^{\leq 1}) \otimes \mcO_{\mbP^1}(-1) \longmigi 0.
\end{align}
Here, observe that 
 $(\mcD_{\mbP^1}^{\leq 2} /\mcD_{\mbP^1}^{\leq 1}) \otimes \mcO_{\mbP^1}(-1) \cong \mcO_{\mbP^1}(3)$ and $\mcD_{\mbP^1}^{\leq 1} \otimes \mcO_{\mbP^1}(-1)$ is an extension of $\mcO_{\mbP^1}(1)$ by $\mcO_{\mbP^1}(-1)$.
Hence, $\mr{Hom}_{\mcO_{\mbP^1}}((\mcD_{\mbP^1}^{\leq 2} /\mcD_{\mbP^1}^{\leq 1}) \otimes \mcO_{\mbP^1}(-1), \mcD_{\mbP^1}^{\leq 1} \otimes \mcO_{\mbP^1}(-1)) =0$, which implies that
there is no splitting of (\ref{E009ju}) but the splitting corresponding to $D_0^\clubsuit$.
This completes the proof of the assertion.
\end{proof}
\vspace{3mm}

\begin{rema} \label{R001}
\leavevmode\\
\ \ \ 
By an argument similar to the argument in the proof of Proposition \ref{P205},
we can prove that both $\mfP \mfC_{\mbP^1, \mbL_0}^n$ and $\mfP \mfC_{\mbP^1, \mbL_0}^{n, \mr{full}}$ (for any $n \geq 2$) are singletons.
In particular, the unique element can be obtained as the image $\xi^{2 \migi n} (D_0^\clubsuit)$ of $D_0^\clubsuit$ via $\xi^{2 \migi n}$ (cf. (\ref{E4289})).
On the open subscheme 
$U := \mr{Spec}(k[u])$ (resp., $T:= \mr{Spec}(k[t])$) as before,
the operator $\xi^{2 \migi n} (D_0^\clubsuit) |_{U} : \mcO_{\mbP^1} (n-1)|_U \migi \mcO_{\mbP^1}(-n-1)|_U$ (resp., $\xi^{2 \migi n} (D_0^\clubsuit) |_{T} : \mcO_{\mbP^1} (n-1)|_T \migi \mcO_{\mbP^1}(-n-1)|_T$) may be expressed as
\begin{align}
f (u) \cdot y^{n-1} \migi \frac{\partial^n f}{\partial u^n} (u) \cdot y^{-n-1} \ \left(\text{resp.,} \ g (t) \cdot x^{n-1} \migi \frac{\partial^n g}{\partial t^n} (t) \cdot x^{-n-1} \right).
\end{align} 
\end{rema}

\vspace{10mm}
\section{The main theorem} \vspace{3mm}

The fourth section is devoted to state and prove the main theorem of the present paper.

\vspace{5mm}
\subsection{Statement of the main theorem} \label{S1g1}
\leavevmode\\ \vspace{-4mm}

Let us fix a smooth curve $X$,  a theta characteristic  $\mbL := (\mcL, \psi_\mcL)$ on $X$.  
The morphism 
\begin{align} \label{W200}
\psi^\A_{\mcL} : \A(\mcL)^\times \migi \A(\Omega_X)^\times
\end{align}
over $X$ between algebraic surfaces  defined by $\psi^\A_{\mcL} (v) = \psi_\mcL (v \otimes v)$  (for each local section $v \in \mcL$)
is  a Galois
 double covering whose Galois group is isomorphic to $\mu_2 = \{ \pm 1 \}$.
 The automorphism of $\A(\mcL)^\times$ determined by   $(-1) \in \mu_2$ is  
  given by $v \mapsto - v$.
Let us write 
\begin{align} \label{Erww23}
\check{\omega}_{\mbL} : =  (\psi_{\mcL}^\A)^*\left(\check{\omega}^{\mr{can}}\right)
\in \Gamma (\A(\mcL)^\times, \bigwedge^2 \Omega_{\A(\mcL)^\times})
\end{align}
(cf. (\ref{W104})), which 
  specifies a symplectic structure on $\A(\mcL)^\times$.
In particular, the pair 
\begin{align}
(\A(\mcL)^\times, \check{\omega}_{\mbL}).
\end{align}
forms a symplectic variety equipped with a $\mu_2$-action.
Then, the main result of the present paper is as follows.

\vspace{3mm}
\bt  \label{T0090}\leavevmode\\
 \ \ \ 
There exists a canonical construction of a $\mu_2$-FC quantization on $(\A (\mcL)^\times, \check{\omega}_\mbL)$ by means of  a Frobenius-projective structure (or equivalently, a dormant indigenous bundle, a dormant $(\mr{SL}_2, \mbL)$-oper,  or a projective connection for $\mbL$ having a full set of solutions)  on $X$. 
The resulting map  of sets
\begin{align}
\bigstar_{X, \mbL} : \mfP \mfS^\mr{F}_X \migi \mfQ^{\mu_2 \text{-}\mr{FC}}_{(\A(\mcL)^\times, \check{\omega}_\mbL)}
\end{align}
is injective and the composite injection
\begin{align}
\bigstar_X : \mfP \mfS^\mr{F}_X \migiincl \mfQ^\mr{FC}_{(\A(\Omega_X)^\times, \check{\omega}^\mr{can})}
\end{align}
of this map and  the natural bijection $\mfQ^{\mu_2 \text{-}\mr{FC}}_{(\A(\mcL)^\times, \check{\omega}_\mbL)} \isom \mfQ^{\mr{FC}}_{(\A(\Omega_X)^\times,  \check{\omega}^\mr{can})}$ (i.e., the inverse of (\ref{W1005})) does not depend on the choice of the theta characteristic $\mbL$.
 \et
\vspace{3mm}

In 
the rest of this section, we shall will  prove the above theorem.

\vspace{5mm}
\subsection{Step I: Local construction} \label{S11}
\leavevmode\\ \vspace{-4mm}

Let 
$\mcS^\heartsuit$
 be a Frobenius-projective structure on $X$, and denote by $\mcE^\spadesuit := (\mcE, \nabla_\mcE, \sigma_\mcE)$ and  $\mcF^\diamondsuit := (\mcF, \nabla_\mcF, \mcN_\mcF, \eta_\mcF)$ the corresponding indigenous bundle and  $(\mr{SL}_2, \mbL)$-oper  respectively.
In this first step, we construct FC  quantizations on various open subschemes of $(\A (\mcL)^\times, \check{\omega}_\mbL)$.
Since $\nabla_\mcF$ has vanishing $p$-curvature, 
$(\mcF, \nabla_\mcF)$ is locally trivial.
More precisely,
there exists a collection
\begin{align} \label{E004}
\{ (U_\alpha, \gamma_\alpha ) \}_{\alpha \in I}
\end{align}
of pairs $(U_\alpha, \gamma_\alpha)$  (indexed by a set $I$), where
$\{ U_\alpha \}_{\alpha \in I}$ is   an open covering  of $X$ and
$\gamma_\alpha$ (for each $\alpha  \in I$) denotes an $\mcO_{U_\alpha}$-linear  isomorphism $\mcF |_{U_\alpha} \isom \mcO_{U_\alpha}^{\oplus 2}$ 
inducing, via taking determinants, the isomorphism $\delta_{\mcF^\diamondsuit} |_{U_\alpha}$ (cf.  (\ref{Effgh3})).
Each  isomorphism $\gamma_\alpha$  induces  an isomorphism 
\begin{align}
 \gamma^\mbA_{\alpha} : \A (\mcF |_{U_\alpha})^\times \isom  \left(\A (\mcO_{U_\alpha}^{\oplus 2})^\times\isom \right) \  U_\alpha \times \mbA^{2 \times}   \end{align}
over $U_\alpha$, where $\mbA^{2 \times} := \mr{Spec}(k[x, y])$ and the isomorphism in the parenthesis is given by $x \mapsto (1, 0)$, $y \mapsto (0, 1)$.
Moreover, 
denote by $\Phi_\alpha$ the composite
\begin{align} \label{Ert52}
\Phi_\alpha : \A (\mcL |_{U_\alpha})^\times \xrightarrow{\eta_\mcF^\A |_{U_\alpha}} \A (\mcF |_{U_\alpha})^\times \xrightarrow{\gamma_{\alpha, \A}}
 U_\alpha \times \mbA^{2 \times} \xrightarrow{\mr{pr}_2} \mbA^{2 \times},
\end{align}
which is compatible with the respective $\mu_2$-actions on $\A (\mcL |_{U_\alpha})^\times$ and $\mbA^{2 \times}$.
The morphism $\A (\Omega_{U_\alpha})^\times \isom \mbA^{2 \times}_{/\mu_2}$ induced  via quotient does not depend on the choice of $\mbL$.

\vspace{3mm}
\ble \label{Lff1}\leavevmode\\
 \ \ \ 
The morphism $\Phi_\alpha$ is \'{e}tale and satisfies the equality  $\Phi_\alpha^*(\omega^{\mr{Weyl}}) = \check{\omega}_\mbL |_{U_\alpha}$. 
 \ele
\begin{proof}
To begin with, we shall consider the case where $X =U_\alpha = \mbP^1 \ (= \mr{Proj}(k [x, y]))$, $\mbL = \mbL_0$, $\mcF^\diamondsuit = \mcF_0^\diamondsuit$, and $\gamma_\alpha$ is the identity of $\mcO_{\mbP^1}^{\oplus 2}$.
Write $\Phi_0$ for the morphism ``$\Phi_\alpha$" in this case.
On the open subscheme $U := \mr{Spec}(k[u])$ of $\mbP^1$,
 where $u = x/y$, 
 we have $\mbA(\mcO_{\mbP^1}(-1))^\times |_{U} = \mr{Spec}(k[u, y, y^{-1}])$, $\mbA^{2 \times}|_U = \mr{Spec}(k[x, y, y^{-1}])$ and $\Phi_0 |_U$ may be given by
  $x \mapsto u \cdot y$ and $y \mapsto y$.
Hence, 
\begin{align} \label{E8999}
\Phi_0^* (\omega^\mr{Weyl})|_U = \Phi^*_0 (dy \wedge dx)= dy \wedge d(u \cdot y) = y \cdot dy \wedge du.
\end{align}
On the other hand, it follows from the definition of $\psi_0$ that
$\psi_0^\mbA : \A (\mcO_{\mbP^1}(-1))^\times \migi \A (\Omega_{\mbP^1})^\times$ is given by assigning $f \mapsto  (y \cdot f)^2\cdot du$ for each $f \in \Gamma (U, \mcO_{\mbP^1}(-1)^\times)$.
If $u^\vee$ denote the dual coordinate of $u$ in $\A (\Omega_{\mbP^1})^\times$,
then 
\begin{align} \label{E8989}
\check{\omega}_{\mbL_0}|_{U} = (\psi^\mbA_{0})^*(\check{\omega}^\mr{can}|_{U}) = (\psi^\mbA_{0})^*\left(\frac{1}{2} \cdot d u^\vee \wedge du\right) = \frac{1}{2} \cdot dy^2 \wedge du = y \cdot dy \wedge du.
\end{align}
By (\ref{E8999}) and (\ref{E8989}), we obtain the desired equality
\begin{align} \label{Ertgbn}
\Phi_0^* (\omega^\mr{Weyl}) = \check{\omega}_{\mbL_0}
\end{align}
 of this case.
(This result will be used in (\ref{Fgghyn}).)

Now, let us go back to our situation.
Denote by 
\begin{align}
\gamma^\mbP_{\alpha} : \mcE |_{U_\alpha} \isom U_\alpha \times \mbP^1
\end{align}
the isomorphism induced from $\gamma^\mbA_{\alpha}$ via projectivization.
Denote by $\phi_\alpha$ the composite
\begin{align}
\phi_\alpha : U_\alpha \xrightarrow{\sigma_\mcE |_{U_\alpha}} \mcE |_{U_\alpha} \xrightarrow{\gamma_{\alpha}^\mbP} U_\alpha \times \mbP^1 \xrightarrow{\mr{pr}_2} \mbP^1,
\end{align}
which is verified to be \'{e}tale since the Kodaira-Spencer map 
associated to $\sigma_\mcE$  is an isomorphism. 
Under the natural identifications 
$U_\alpha \times \mbA^{2 \times} \isom \phi_\alpha^*(\A (\mcO_{\mbP^1}^{\oplus 2})^\times)$ 
and
 $U_\alpha \times \mbP^1 \isom \phi_\alpha^*(\mcE_0)$ (where $\phi_\alpha^*(-)$ denotes base-change by $\phi_\alpha$), we obtain a commutative diagram
\begin{align}
\vcenter{\xymatrix{ &  & & \A (\mcF |_{U_\alpha})^\times \ar[dd]^{\pi_{\mcF|_{U_\alpha}}} \ar[dl]^\sim_{\gamma^\mbA_{\alpha}} \\
& & \phi_\alpha^*(\A (\mcO_{\mbP^1}^{\oplus 2})^\times) \ar[dd]^(.35){\phi_\alpha^*(\pi_{\mcO_{\mbP^1}^{\oplus 2}})}& \\
& U_\alpha \ar[dl]^\sim_{\mr{id}_{U_\alpha}} \ar[rr]^{\hspace{-25mm}\sigma_\mcE |_{U_\alpha}}|(.45)\hole& & \mcE |_{U_\alpha} \ar[ld]_\sim^{\gamma_{\alpha}^\mbP} \\
U_\alpha \ar[rr]_{\phi_\alpha^*(\sigma_0)}& & \phi_\alpha^*(\mcE_0).& 
}}
\end{align}
Since   (\ref{Ddd01}) is cartesian, the above diagram induces an isomorphism (of $\mbG_m$-torsors)
\begin{align} \label{E3321}
\underline{\gamma}_{\alpha}^\mbA : \A (\mcL |_{U_\alpha})^\times \isom \phi_\alpha^*(\mbA (\mcO_{\mbP^1}(-1))^\times)
\end{align}
over $U_\alpha$ such that the following diagram is commutative:
\begin{align} \label{Ertyui}
\vcenter{\xymatrix{ & \A (\mcL |_{U_\alpha})^\times \ar[rr]^{\eta^\mbA_\mcF |_{U_\alpha}} \ar[dl]_{\underline{\gamma}_{\alpha}^\mbA}^\sim \ar[dd]|(.50)\hole& & \A (\mcF |_{U_\alpha})^\times \ar[dd]^{\pi_{\mcF|_{U_\alpha}}} \ar[dl]_\sim^{\gamma_{\alpha}^\mbA} \\
\phi^*_\alpha (\A (\mcO_{\mbP^1}(-1))^\times) \ar[rr]^{\hspace{20mm}\phi^*_\alpha (\eta^\mbA_{0})} \ar[dd] & & \phi_\alpha^*(\A (\mcO_{\mbP^1}^{\oplus 2})^\times) \ar[dd]^(.35){\phi_\alpha^*(\pi_{\mcO_{\mbP^1}^{\oplus 2}})}& \\
& U_\alpha \ar[dl]^-\sim_-{\mr{id}_{U_\alpha}} \ar[rr]^{\hspace{-20mm}\sigma_\mcE |_{U_\alpha}}|(.50)\hole& & \mcE |_{U_\alpha} \ar[ld]_\sim^{\gamma_{\alpha}^\mbP} \\
U_\alpha \ar[rr]_{\phi_\alpha^*(\sigma_0)}& & \phi_\alpha^*(\mcE_0).& 
}}
\end{align}
Moreover, the morphism (\ref{E3321}) induces an $\mcO_{U_\alpha}$-linear isomorphism 
\begin{align} \label{E4291}
\underline{\gamma}_\alpha : \mcL|_{U_\alpha} \isom \phi_\alpha^*(\mcO_{\mbP^1}(-1))
\end{align}
 fitting into  the following isomorphism of short exact sequences:
\begin{align}
\vcenter{\xymatrix{
0 \ar[r]& \mcL |_{U_\alpha} \ar[r]^{\eta_\mcF} \ar[d]_\wr^{\underline{\gamma}_\alpha}& \mcF |_{U_\alpha} \ar[d]_\wr^{\gamma_\alpha} \ar[r] & \mcL^\vee |_{U_\alpha} \ar[r] \ar[d]_\wr^{(\underline{\gamma}_\alpha^{\vee})^{-1}}& 0 \\
 0 \ar[r] & \phi^*_\alpha (\mcO_{\mbP^1}(-1)) \ar[r] & \phi_\alpha (\mcO_{\mbP^1}^{\oplus 2}) \ar[r] & \phi^*_\alpha (\mcO_{\mbP^1}(-1)) \ar[r]& 0,
}}
\end{align}
where the upper right-hand horizontal arrow $\mcF |_{U_\alpha} \migi \mcL^\vee |_{U_\alpha}$ arises from $\delta_{\mcF^\diamondsuit} : \mr{det}(\mcF) \isom \mcO_X$ (cf. (\ref{Effgh3})). 
Since $\gamma_\alpha$ is, moreover,  compatible with the connections $\nabla_\mcF$ and $d^{\oplus 2}$, the respective Kodaira-Spencer maps 
 give rise to
 a commutative square diagram of $\mcO_{U_\alpha}$-modules
\begin{align}
\vcenter{\xymatrix@C=40pt{
\mcL^{\otimes 2} |_{U_\alpha} \ar[r]_\sim^{\psi_{\mcL}|_{U_\alpha}} \ar[d]^\wr_{\underline{\gamma}_\alpha^{\otimes 2}} & \Omega_{U_\alpha} \ar[d]_\wr^{\phi_\alpha}  \\
\phi_\alpha^*(\mcO_{\mbP^1}(-1)^{\otimes 2}) \ar[r]^-\sim_-{\phi_\alpha^*(\psi_{0})} &\phi_\alpha^*(\Omega_{\mbP^1}), 
}}
\end{align}
as well as a commutative diagram of $U_\alpha$-schemes
\begin{align} \label{E687ff}
\vcenter{\xymatrix@C=44pt{
\A (\mcL |_{U_\alpha})^\times  \ar[d]^{\wr}_{\underline{\gamma}_\alpha^\mbA} \ar[r]^{\psi_\mcL^\mbA |_{U_\alpha}}& \A  (\Omega_{U_\alpha})^\times  \ar[d]_\wr^{\phi_\alpha^{\mbA}} \\
\phi_\alpha^*(\A (\mcO_{\mbP^1}(-1))^\times) \ar[r]_-{\phi_\alpha^*(\psi_{\mbL_0}^\mbA)}& \phi_\alpha^*(\mbA (\Omega_{\mbP^1})). 
}}
\end{align}
Hence, the following sequence of equalities holds:
\begin{align} \label{Fgghyn}
\Phi_\alpha^*(\omega^\mr{Weyl}) &\stackrel{(\ref{Ert52})}{=} (\gamma_{\alpha}^\mbA \circ \eta_\mcF^\mbA |_{U_\alpha})^*(\mr{pr}_2^*(\omega^{\mr{Weyl}})) \\
&\stackrel{(\ref{Ertyui})}{=} (\phi_\alpha^*(\eta_0^\mbA) \circ \underline{\gamma}_{\alpha}^\mbA)^*(\mr{pr}_2^*(\omega^\mr{Weyl})) \notag \\
&= \underline{\gamma}_{\alpha}^{\mbA*}(\Phi_0^*(\omega^\mr{Weyl})) \notag \\
&\stackrel{(\ref{Ertgbn})}{=} \underline{\gamma}_{\alpha}^{\mbA*} (\check{\omega}_{\mbL_0}) \notag \\
&\stackrel{(\ref{Erww23})}{=}  \underline{\gamma}_{\alpha}^{\mbA*} \left(\frac{1}{2} \cdot (\psi_0^\mbA)^*(\omega^\mr{can}_{\mbP^1})\right) \notag \\
&\stackrel{(\ref{E687ff})}{=} \frac{1}{2} \cdot (\phi_\alpha^{\mbA*}\circ\psi_\mcL^\mbA)^*(\omega_{\mbP^1}^\mr{can}) \notag \\
& = \frac{1}{2} \cdot \psi_{\mcL}^{\mbA*}(\omega_{U_\alpha}^\mr{can}) \notag \\
& = \check{\omega}_\mbL |_{U_\alpha}. \notag
\end{align}
This completes the proof of the latter assertion.

The former assertion, i.e., the \'{e}taleness of $\Phi_\alpha$ follows from the
\'{e}taleness of $\phi_\alpha$ and the fact that the square diagram
\begin{align}
\vcenter{\xymatrix{
\mbA(\mcL |_{U_\alpha})^\times \ar[r]^-{\Phi_\alpha} \ar[d] & \mbA^{2 \times} \ar[d]\\
U_\alpha \ar[r]_{\phi_\alpha} & \mbP^1
}}
\end{align}
is cartesian, where the vertical arrows are the natural projections.
This completes the proof of the lemma.
\end{proof}
\vspace{3mm}

By  the above lemma and the discussion in \S\,\ref{S10} applied to $\Phi_\alpha$,
the pull-back of the $\mu_2$-FC quantization $\mcW^2_k$ (cf. (\ref{Efgklo})) on $(\mbA^{2\times}, \omega^\mr{Weyl})$ specifies 
 a $\mu_2$-FC quantization
\begin{align} \label{W201}
\Phi^*_\alpha(\mcW^2_k)
\end{align}
on the symplectic variety $(\A (\mcL |_{U_\alpha})^\times, \check{\omega}_\mbL |_{U_\alpha})$.
The (isomorphism class of the) FC quantization on $(\A (\Omega_{U_\alpha}), \check{\omega}^\mr{can})$ 
corresponding to 
$\Phi^*_\alpha(\mcW^2_k)$ via (\ref{W1005}) does not depend on the choice of $\mbL$.

\vspace{5mm}
\subsection{Step II: Global construction} \label{S11}
\leavevmode\\ \vspace{-4mm}

In this second step, we glue together the locally defined quantizations constructed  above to obtain an FC quantization on the entire space $X$, as follows.
After possibly replacing $\{U_\alpha \}_\alpha$ with its refinement, we can assume that each $U_\alpha$ is affine.
Let us take a pair  $(\alpha, \beta) \in I \times I$ with $U_{\alpha\beta} := U_\alpha \cap U_\beta \neq \emptyset$.
Since ($X$ is separated, which implies that) $U_{\alpha \beta}$ is affine, we can write $U_{\alpha \beta} = \mr{Spec}(R_{\alpha \beta})$ for some $k$-algebra  $R_{\alpha \beta}$.
In what follows, we shall use the notation $(-)^{(1)}$ to denote the base-change of objects via $F_k$.  
In particular, we obtain the $k$-algebra $R^{(1)}_{\alpha \beta}$ equipped with a $k$-algebra (injective) homomorphism 
$R^{(1)}_{\alpha \beta} \migi R_{\alpha \beta}$.
Also, write $R_{\mcL, \alpha \beta} := \Gamma (\A (\mcL |_{U_{\alpha \beta}})^\times, \mcO_{\A (\mcL |_{U_{\alpha \beta}})^\times})$.
If $\mcI_\alpha$ denotes  the ideal of $R_{\alpha \beta}^{}[x, y]$ determined by the closed immersion
\begin{align}
(\gamma_{\alpha}^\mbA\circ \eta_\mcF^\mbA|_{U_\alpha}) |_{U_{\alpha \beta}}: 
\A (\mcL |_{U_{\alpha \beta}})^\times \ \left(= \mr{Spec}(R_{\mcL, \alpha \beta}) \right) \migiincl \mbA^{2 \times} \times U_{\alpha \beta} \ \left(= \mr{Spec} (R_{\alpha \beta}[x, y]) \right).
\end{align}
Hence, $R_{\alpha \beta}[x, y]/\mcI_\alpha \cong R_{\mcL, \alpha \beta}$ and
 we have a natural isomorphism of $R^{(1)}_{\alpha \beta}[x^p, y^p]/\mcI_\alpha^{(1)}$ ($= \Gamma ((\A (\mcL |_{U_{\alpha \beta}})^\times)^{(1)}, \mcO_{(\A (\mcL |_{U_{\alpha \beta}})^\times)^{(1)}})$)-algebras
\begin{align} \label{W210}
\Gamma (\A (\mcL |_{U_{\alpha \beta}})^\times, \Phi_\alpha^*(\mcW_k^2))  \cong W^2_{R_{\alpha \beta}} \otimes_{R_{\alpha \beta}^{(1)}[x^p, y^p]} (R^{(1)}_{\alpha \beta}[x^p, y^p]/\mcI_\alpha^{(1)}).
\end{align}

Next,
let $\gamma_{\alpha \beta}$  be the automorphism of $R_{\alpha \beta}[x, y]$
corresponding to  
  $(\gamma_\alpha |_{U_{\alpha \beta}}) \circ (\gamma_\beta |_{U_{\alpha \beta}})^{-1} \in \mr{SL}_2 (R_{\alpha \beta}) \ \left(= \mr{Sp}(R_{\alpha \beta}) \right)$.
This automorphism restricts to an automorphism of $R_{\alpha \beta}[x^p, y^p]$, and hence, we have the following diagram
\begin{align} \label{E02jj5}
\xymatrix@R=6pt{
&& R_{\mcL, \alpha \beta} &&
\\
&R_{\alpha \beta}[x, y] \ar[ru]&& R_{\alpha \beta}[x, y] \ar[ll]^{\hspace{20mm}\gamma_{\alpha \beta}}_{\hspace{19mm}\sim} \ar[lu]&
\\
&& R_{\mcL, \alpha \beta}^{(1)} \ar[uu]|(.51)\hole &&
\\
&R_{\alpha \beta}^{(1)}[x^p, y^p]  \ar[uu] \ar[ru] && R_{\alpha \beta}^{(1)}[x^p, y^p]. \ar[ll]^{\gamma_{\alpha \beta}}_\sim  \ar[uu]\ar[lu]&
}
\end{align}
The bottom triangle in (\ref{E02jj5}) turns out to be commutative since the other various small diagrams are commutative.
This implies that 
$\gamma_{\alpha \beta} (\mcI_\alpha^{(1)}) = \mcI_\beta^{(1)}$.
Hence, the automorphism
of the $R_{\alpha \beta}^{(1)}$-algebras $W^2_{R_{\alpha \beta}}$ given by 
$\gamma_{\alpha \beta}$ via (\ref{E4522})
 induces an isomorphism
\begin{align}
W^2_{R_{\alpha \beta}} \otimes_{R_{\alpha \beta}^{(1)}[x^p, y^p]} (R_{\alpha \beta}^{(1)}[x^p, y^p]/\mcI_\beta^{(1)}) \isom W^2_{R_{\alpha \beta}} \otimes_{R_{\alpha \beta}^{(1)}[x^p, y^p]} (R_{\alpha \beta}^{(1)}[x^p, y^p]/\mcI_\alpha^{(1)})
\end{align}
By passing to (\ref{W210}), we obtain an isomorphism
\begin{align}
\Phi_{\alpha \beta} : \Phi_\beta^*(\mcW_k^2) |_{\A (\mcL |_{U_{\alpha \beta}})^\times} \isom \Phi_\alpha^*(\mcW_k^2)|_{\A (\mcL |_{U_{\alpha \beta}})^\times}
\end{align}
of sheaves of $k[[\hslash]]$-algebras on $\A (\mcL |_{U_{\alpha \beta}})^\times$.
This isomorphism is verified to be compatible with the $\mu_2$-actions,  and the induced isomorphism between the respective quotient FC quantizations does not depend on the choice of $\mbL$. 

By means of the isomorphisms $\Phi_{\alpha \beta}$ (for various $\alpha$, $\beta$), 
the sheaves $\Phi_\alpha^*(\mcW_k^2)$
 may be glued together to obtain  a sheaf 
 \begin{align}
 \mcW_{\mcS^\heartsuit},
 \end{align}
which  forms a $\mu_2$-FC quantization on $(\A (\mcL)^\times, \check{\omega}_\mbL)$.
 The isomorphism class of this quantization does not depend on the choice of  $\{ (U_\alpha, \gamma_\alpha) \}_{\alpha \in I}$, and moreover, the isomorphism class of its quotient by the $\mu_2$-action does not depend on the choice of $\mbL$.
Thus, we have obtained a well-defined map
 \begin{align}
\bigstar_{X, \mbL}  :  \mfP \mfS^F_X \migi \mfQ^{\mu_2 \text{-}\mr{FC}}_{(\A (\mcL)^\times, \check{\omega}_\mbL)}
 \end{align}
such that the composite 
\begin{align}
\bigstar_X  :  \mfP \mfS^F_X \migi \mfQ^{\mr{FC}}_{(\A (\Omega_X)^\times, \check{\omega}^\mr{can})}
\end{align}
of this map and  the  bijection $\mfQ_{(\A (\mcL)^\times, \check{\omega}_\mbL)}^{\mu_2 \text{-}\mr{FC}}\isom \mfQ^\mr{FC}_{(\A (\Omega_X)^\times,  \check{\omega}^\mr{can})}$  does not depend on the choice of $\mbL$.

\vspace{5mm}
\subsection{Step III: Injectivity} \label{S62}
\leavevmode\\ \vspace{-4mm}

The remaining portion of the proof is the injectivity of $\bigstar_{X, \mbL}$.
Here, let $(S, \omega_S)$ be a symplectic variety and $\mcO_S^\hslash$ an FC quantization on $(S, \omega_S)$.
Given two local sections $a, b \in \mcO_S$, we express $a * b$ as 
\begin{align}
a * b = \delta_b^0(a) + \delta_b^1 (a) \cdot  \hslash + \delta_b^2 (a) \cdot  \hslash^2 + \cdots  \in \mcO_{S} [[\hslash ]]
\end{align}
for some local sections $\delta_b^i (a) \in \mcO_S$ ($i=0,1,2, \cdots$).
For each $i$, the assignment $b \mapsto \delta_b^i (a)$ defines a (locally defined) 
$k$-linear endomorphism  $\delta_b^i \in \mcE nd_k (\mcO_{S})$.
The assignment $b \mapsto \delta_b^i$ gives an $\mcO_{S}$-linear morphism
\begin{align}
\delta^i : \mcO_{S} \migi \mcE nd_k (\mcO_{S}),
\end{align}
where we equip $\mcE nd_k (\mcO_{S})$ with a structure of $\mcO_S$-module given by multiplication on the left.

\vspace{3mm}
\ble \label{L200} \leavevmode\\
 \ \ \ 
If the triple $(S, \omega_S, \mcO_S^\hslash)$ is taken to be $(\A (\mcL)^\times, \check{\omega}_\mbL, \mcW_{\mcS^\heartsuit})$ (as discussed in the previous subsection),
then
for each $i < p$  the image of $\delta^i$ is contained in $\mcD^{\leq i}_{\A (\mcL)^\times}$.
 \ele
\begin{proof}
By the local nature of the assertion, it suffices to prove this lemma with $\A (\mcL)^\times$ replaced with each open subscheme $\A (\mcL |_{U_\alpha})^\times$.
Moreover, since
$\Phi_\alpha^*({^l \mcD}^{\leq i}_{\A^{2\times}})\cong {^l \mcD}^{\leq i}_{\A (\mcL|_{U_\alpha})^\times}$,
it suffices to consider the case where $(S, \omega_S, \mcO_S^\hslash)$ is taken to be $(\mbA^{2 \times}, \omega^\mr{Weyl}, \mcW_k^2)$.
 Then, it follows from the definition of the multiplication in $W_k^2$ that for each $a, b \in k[x, y]$, we have $\delta^i_b(a) = \frac{1}{i!} \cdot \frac{\partial^i a}{\partial y^i}\cdot  \frac{\partial^i b}{\partial x^i}$.
 This completes the proof of the lemma.
\end{proof}
\vspace{3mm}

We shall finish the proof of the main theorem.
Since $\pi_*( \mcO_{\A (\mcL)^\times})$ is naturally identified with $\bigoplus_{j \in \mbZ}\mcL^{\otimes j}$, we have the inclusion
\begin{align} \label{Errt}
\mcL^{\otimes (-i+1)} \migiincl
 \pi_*( \mcO_{\A (\mcL)^\times})
\end{align}
into the $(-i+1)$-st factor.
Next, if $\pi$ denotes the natural projection $\mbA (\mcL)^\times \migi X$, then
 the kernel of  the surjection $\mcT_{\A (\mcL)^\times} \migisurj\pi^*(\mcT_{X})$ 
obtained by differentiating $\pi$ is isomorphic to $\pi^*(\mcL)$.
The resulting injection $\pi^*(\mcL) \migiincl \mcT_{\A (\mcL)^\times}$ induces an injection
\begin{align} \label{Ertygfd}
\pi^*(\mcL^{\otimes i}) \migiincl \mcT_{\A (\mcL)^\times}^{\otimes i}.
\end{align}
The projection formula gives the composite isomorphism
\begin{align} \label{F9090}
\pi_*(\pi^*(\mcL^{\otimes i})) \isom \pi_*(\mcO_{\A (\mcL)^\times})\otimes \mcL^{\otimes i} \isom \left(\bigoplus_{j \in \mbZ} \mcL^{\otimes j}\right) \otimes \mcL^{\otimes i} = \bigoplus_{j \in \mbZ} \mcL^{\otimes (j+i)}.
\end{align}
Hence, we have an injection
\begin{align} \label{Eghj}
\mcL^{\otimes (i+1)} \migiincl \bigoplus_{j \in \mbZ} \mcL^{\otimes (j+i)} \xrightarrow{(\ref{F9090})} \pi_*(\pi^*(\mcL^{\otimes i})) \xrightarrow{(\ref{Ertygfd})}\pi_*(\mcT_{\A (\mcL)^\times}^{\otimes i}), 
\end{align}
where the first arrow is the inclusion into the first factor.

\vspace{3mm}
\bpr \label{L203} \leavevmode\\
 \ \ \ 
Let $i$ be an integer with $2 \leq i < p$ and let $\delta^i : \mcO_{\mbA (\mcL)^\times} \migi \mcD_{\mbA (\mcL)^\times}^{\leq i}$ be the morphism resulting from Lemma \ref{L200}.
Denote by $D^\clubsuit$  the projective connection for $\mbL$ corresponding to $\mcS^\heartsuit$ via the upper horizontal bijections in the diagram (\ref{E46090}).
Then, the following  diagram is commutative:
\begin{align} \label{Efgmaob}
\vcenter{\xymatrix@C=46pt{
\mcL^{\otimes (-i +1)} \ar[rr]^{\frac{1}{i!} \cdot \xi^{2 \migi i} (D^\clubsuit)} \ar[d]_{(\ref{Errt})} & & \mcL^{\otimes (i+1)} \ar[d]^{(\ref{Eghj})}
\\
\pi_*(\mcO_{\A (\mcL)^\times})  \ar[r]_{\pi_*( \delta^i)}  & \pi_* (\mcD^{\leq i}_{\mbA (\mcL)^\times}) \ar[r]_{\pi_*(\Sigma)} &\pi_*(\mcT^{\otimes i}_{\A (\mcL)^\times}).
}}
\end{align}
In particular, (by considering the case of $i=2$), the map $\bigstar_X$  is injective.

 \epr
\begin{proof}
First, we shall consider the former assertion.
Just as in the proof of Lemma \ref{L200},
it suffices to consider the case where the triple  $(X, \mbL, \mcS^\heartsuit)$ is taken to be $(\mbP^1, \mbL_0, \mcS_0^\heartsuit)$.
Let us  identify $\mbA (\mcO_{\mbP^1}(-1))^\times$ with $\mbA^{2 \times}$ via $\Phi_0$ (cf. the proof in Lemma \ref{Lff1} for the definition of $\Phi_0$).
Over the open subsheme $U:= \mr{Spec}(k[u])$ (where $u:= x/y$) of $\mbP^1$,
the composite $\Sigma \circ \delta^i$ sends each element $f (u) \cdot y^{i-1} \in \Gamma (U, \mcO_{\mbP^1}(i-1))$ to
\begin{align} \label{Efghhy}
\frac{1}{i!} \cdot \frac{\partial^i}{\partial x^i}(f(u) \cdot y^{i-1}) \cdot \left(\frac{\partial}{\partial y} \right)^{\otimes i} = \frac{1}{y \cdot i!} \cdot \frac{\partial^if}{\partial u^i}(u) \cdot \left(\frac{\partial}{\partial y} \right)^{\otimes i} \in \Gamma (\pi^{-1}(U), \mcT_{\mbA^{2\times}}^{\otimes i}).
\end{align}
On the other hand, since the injection $\mcO_{\mbP^1}(-i-1) \migi \pi_*(\mcT^{\otimes i}_{\mbA (\mcO_{\mbP^1}(-1))^\times})$ of (\ref{Eghj}) is given by $\frac{1}{y^{i+1}} \mapsto \frac{1}{y} \cdot \big(\frac{\partial}{\partial y}\big)^{\otimes i}$, the image of 
$\frac{1}{i!} \cdot \xi^{2\migi u}(D_0^\clubsuit)(f(u) \cdot y^{i-1})$ via this injection coincides with $\frac{1}{y \cdot i!}\cdot \frac{\partial^i f}{\partial u^i} (u) \cdot \big(\frac{\partial}{\partial y}\big)^{\otimes i}$ (cf. Remark \ref{R001}), which is identical to (\ref{Efghhy}).
This implies the commutativity of  (\ref{Efgmaob}), as desired.
The latter assertion follows from the former assertion.
\end{proof}
\vspace{3mm}

According to the above results, we complete the proof of the main theorem.

\vspace{10mm}
\section{A higher-dimensional variant of the main theorem}  \label{S010}\vspace{3mm}

In this final section, we shall prove (cf. Theorem \ref{T009h1gh} described later) a higher dimensional variant of our main theorem, which may be thought of as a positive characteristic analogue of ~\cite{Bi}, Proposition 4.3.

\vspace{5mm}
\subsection{Frobenius-$\mr{Sp}$ structures} \label{Sp62}
\leavevmode\\ \vspace{-4mm}

Let $n$ be an integer with $n >1$ and 
$S$  a smooth variety of dimension $2n-1$. 

\vspace{3mm}
\bde
 \label{D0799} \leavevmode\\
 \ \ \ 
 A {\bf Frobenius-$\mr{Sp}$ structure (of level $1$)} on $S$ is a triple
 \begin{align}
 \mcS^{\heartsuit \diamondsuit} := (\mcS^\heartsuit, \mcS^\diamondsuit, \kappa)
 \end{align}
 consisting of a Frobenius-projective structure (of level $1$) $\mcS^\heartsuit$ on $S$,  an $(\mr{Sp}_{2n})_S^\mr{F}$-torsor $\mcS^\diamondsuit$ on $S$, and an isomorphism $\kappa : \mcS^\diamondsuit \times^{\mr{Sp}_{2n}} \mr{PGL}_{2n} \isom \mcS^\heartsuit$, where 
 $\mcS^\diamondsuit \times^{\mr{Sp}_{2n}} \mr{PGL}_{2n}$ denotes 
   the $(\mr{PGL}_{2n})_S^\mr{F}$-torsor induced from $\mcS^\diamondsuit$ via   change of structure group by  the composite $\mr{Sp}_{2n} \migiincl \mr{GL}_{2n} \migisurj \mr{PGL}_{2n}$.
  \ede
\vspace{3mm}

In what follows, we shall observe that each Frobenius-$\mr{Sp}$ structure   gives rise to a theta characteristic  and an FC quantization on a certain symplectic variety.

First, let us consider a procedure for constructing  a theta characteristic by means of  a Frobenius-$\mr{Sp}$ structure.
Let $\mcS^{\heartsuit \diamondsuit} := (\mcS^\heartsuit, \mcS^\diamondsuit, \kappa)$ be a Frobenius-$\mr{Sp}$ structure on $S$.
The $(\mr{SL}_{2n})^\mr{F}_S$-torsor corresponding to $\mcS^\diamondsuit$ via change of structure group by the natural inclusion $\mr{Sp}_{2n} \migiincl \mr{SL}_{2n}$ determines a vector bundle on $S^{(1)}$ with trivial determinant.
If $\mcF_\mcS$ denotes its pull-back via $F_{S/k}$, then $\mcF_\mcS$ admits a canonical connection $\nabla_\mcS$ with vanishing $p$-curvature (cf. ~\cite{Kat}, \S\,5, Theorem 5.1). 
In particular, we have $(\mr{det}(\mcF_\mcS), \nabla_{\mr{det}(\mcF_\mcS)}) \cong (\mcO_S, d)$.
Denote by $\mcE_\mcS$ the $\mbP^{2n-1}$-bundle on $S$ defined to be the pull-back via $F_{S/k}$ of the $\mbP^{2n-1}$-bundle on $S^{(1)}$ corresponding to  $\mcS^\heartsuit$.
The local sections $U \migi U \times \mbP^n$ (for various open subschemes $U$ of $S$) classified by sections in $\mcS^\heartsuit$ may be glued together to obtain a well-defined global section 
 $\sigma_\mcS : S \migi \mcE_\mcS$. 
Then, there exists a unique line subbundle $\mcL_\mcS$ of $\mcF_\mcS$ 
 such that the square diagram
\begin{align}
\vcenter{\xymatrix{
\A (\mcL_\mcS)^\times \ar[r] \ar[d] & \A (\mcF_\mcS)^\times\ar[d]
\\
S \ar[r]_{\sigma_\mcS} & \mcE_\mcS
}}
\end{align}
is commutative and cartesian, where the upper horizontal arrow arises from the inclusion $\mcL_\mcF \migiincl \mcF_\mcS$ and  the vertical arrows denote the natural projections.
By the definitions of $\sigma_\mcS$ and $\mcL_\mcS$, 
the $\mcO_{S}$-linear morphism 
\begin{align}
\Omega_S^\vee \otimes \mcL_\mcS \isom \mcF_\mcS / \mcL_\mcS
\end{align}
 induced naturally by the composite
\begin{align}
\mcL_\mcS \migiincl \mcF_\mcS \xrightarrow{\nabla_\mcS} \Omega_S \otimes \mcF_\mcS \migisurj \Omega_S \otimes (\mcF_\mcS/ \mcL_\mcS)
\end{align}
turns out to be an isomorphism.
This isomorphism yields, via taking determinants, an isomorphism
\begin{align}
\omega_S^\vee \otimes \mcL_\mcS^{\otimes (2n-1)} \ \left(=\mr{det}(\Omega_S^\vee \otimes \mcL_\mcS) \right) \isom \left(\mr{det}(\mcF_\mcS/\mcL_\mcS) = \right) \mcL_\mcS^\vee,
\end{align}
or equivalently, an isomorphism
\begin{align}
\psi_\mcS : \mcL_\mcS^{\otimes 2n} \isom \omega_S.
\end{align}
Thus, we obtain a theta characteristic  
\begin{align}
\mbL := (\mcL_\mcS, \psi_\mcS)
\end{align}
on  $S$.

\vspace{3mm}
\begin{rema} \label{R156}
\leavevmode\\
\ \ \
Let us mention the case  where
 $n=1$, i.e., $S$ is a smooth curve.
 As discussed  above, each Frobenius-$\mr{Sp}$ structure 
  yields a theta characteristic.
 Conversely, suppose that we are given a theta characteristic $\mbL$ and a projective structure $\mcS^\heartsuit$ on $S$.
 Denote by $\mcF^\diamondsuit := (\mcF, \nabla_\mcF, \mcN_\mcF, \eta_\mcF)$  the $(\mr{SL}_2, \mbL)$-oper  corresponding to $\mcS^\heartsuit$.
 Since it has vanishing $p$-curvature, there exists a unique $(\mr{SL}_2)_S^\mr{F}$-torsor $\mcS^\diamondsuit$ such that the induced $\mr{SL}_2$-torsor on $S$ equipped with a canonical connection 
  is isomorphic to $(\mcF, \nabla_\mcF)$.
 Then,   
 (since $\mr{SL}_2 = \mr{Sp}_{2}$) the triple consisting of $\mcS^\heartsuit$, $\mcS^\diamondsuit$, and the natural isomorphism $\kappa : \mcS^\diamondsuit \times^{\mr{SL}_2} \mr{PGL}_2 \isom \mcS^\heartsuit$, specifies a Frobenius-$\mr{Sp}$ structure on $S$.
 According to this construction,  
 giving a Frobenius-$\mr{Sp}$ structure on a smooth curve $S$ is equivalent to giving a pair of a theta characteristic $\mbL$ and a Frobenius-projective structure on $S$.
\end{rema}

\vspace{5mm}
\subsection{An FC quantizations arising from a Frobenius-$\mr{Sp}$ structure} \label{Sj162}
\leavevmode\\ \vspace{-4mm}

Next, let us construct a certain symplectic variety and an FC quantization on it.
We shall keep the above notation.
By the definition of a Frobenius-$\mr{Sp}$ structure, 
there exists a collection 
\begin{align}
\{ (U_\alpha, \gamma_\alpha) \}_{\alpha \in I},
\end{align}
of pairs $(U_\alpha, \gamma_\alpha)$ (indexed by a set $I$),
 where $\{ U_\alpha \}_{\alpha \in I}$ is an open covering of $S$ and $\gamma_\alpha$ (for each $\alpha \in I$) denotes an $\mcO_{U_\alpha}$-linear isomorphism $\mcF_\mcS |_{U_\alpha} \isom \mcO_{U_\alpha}^{\oplus 2n}$ inducing, via taking determinants, the fixed isomorphism $\mr{det}(\mcF_\mcS) \isom \mcO_S$ (restricted to $U_\alpha$).
 Moreover, we can assume that for any pair $(\alpha, \beta) \in I \times I$ with $U_{\alpha \beta} := U_\alpha \cap U_\beta \neq \emptyset$,
  the automorphism $\gamma_{\alpha \beta} := (\gamma_\alpha |_{U_{\alpha \beta}})\circ (\gamma_\beta |_{U_{\alpha\beta}})^{-1}$ of $\mcO_{U_{\alpha \beta}}^{\oplus 2n}$ corresponds to a $U_{\alpha \beta}$-rational point of $\mr{Sp}_{2n} \ \left(\subseteq \mr{SL}_{2n} \right)$.
Let 
\begin{align}
\gamma_\alpha^\A : \A (\mcF_\mcS |_{U_\alpha})^\times
\isom \left( \A (\mcO_{U_\alpha}^{\otimes 2n}) \isom \right) U_\alpha \times \mbA^{2n \times} 
\end{align}
be the isomorphism induced by $\gamma_\alpha$, and let 
$\gamma_\alpha^\mbP$ be the isomorphism $\mcE_\mcS |_{U_\alpha} \isom U_\alpha \times \mbP^{2n-1}$ obtained from $\gamma_\alpha^\mbA$ via projectivization.
Then, we have
two composites
\begin{align}
\Phi_\alpha  &: \mbA (\mcL_\mcS |_{U_\alpha})^\times \xrightarrow{\mr{incl.}} \mbA (\mcF_\mcS |_{U_\alpha})^\times \xrightarrow{\gamma_\alpha^\mbA}   U_\alpha \times \mbA^{2n \times}  \xrightarrow{\mr{pr}_2} \mbA^{2n \times}, \\
\phi_\alpha &: U_\alpha \xrightarrow{\sigma_\mcS |_{U_\alpha}} \mcE_\mcS |_{U_\alpha} \xrightarrow{\gamma_\alpha^\mbP} U_\alpha \times \mbP^{2n -1}  \xrightarrow{\mr{pr}_2} \mbP^{2n-1}. \notag
\end{align}
Since $\phi_\alpha$ is \'{e}tale and the square diagram 
\begin{align}
\vcenter{\xymatrix{
\mbA (\mcL_\mcS |_{U_\alpha})^\times \ar[r]^-{\Phi_\alpha}\ar[d] & \mbA^{2n \times} \ar[d] \\
U_\alpha \ar[r]_-{\phi_\alpha}& \mbP^{2n-1}
}}
\end{align}
(where the vertical arrows denote the natural projections) is commutative and cartesian, 
$\Phi_\alpha$ turns out to be \'{e}tale.
The pull-back  $\Phi_\alpha^* (\omega^\mr{Weyl})$ of  $\omega^\mr{Weyl}$ via $\Phi_\alpha$ speficies
a symplectic structure  on $\mbA (\mcL_\mcS |_{U_\alpha})^\times$.
  For each pair $(\alpha, \beta) \in I \times I$ with $U_{\alpha \beta} \neq \emptyset$, we shall denote by $\gamma_{\alpha \beta}^\mbA$  the automorphism of $U_{\alpha\beta} \times \mbA^{2n \times}$ corresponding to $\gamma_{\alpha \beta}$.
 Since $\gamma_{\alpha \beta} \in \mr{Sp}_2(U_{\alpha \beta})$, the equality  $\gamma_{\alpha \beta}^{\mbA *}(\mr{pr}_2^*(\omega^\mr{Weyl}))= \mr{pr}_2^*(\omega^\mr{Weyl})$ holds,
     which implies that $\Phi_\alpha^*(\omega^\mr{Weyl})|_{U_{\alpha \beta}} = \Phi_\beta^*(\omega^\mr{Weyl})|_{U_{\alpha \beta}}$.
  Thus, the $\Phi_\alpha^*(\omega^\mr{Weyl})$'s may be glued together to obtain a symplectic structure
  $\omega_\mcS$
  on $\mbA (\mcL_\mcS)^\times$.
  In particular, we have a symplectic variety
  \begin{align}
  (\mbA (\mcL_\mcS)^\times, \omega_\mcS).
  \end{align}

 Moreover, for each $\alpha \in I$,
the pull-back $\Phi_\alpha^*(\mcW_k^{2n})$ of $\mcW_k^{2n}$ via $\Phi_\alpha$ specifies 
  an FC quantization
 on $(\mbA (\mcL_\mcS |_{U_\alpha})^\times, \Phi_\alpha^*(\omega^\mr{Weyl}))$.
It follows from an argument similar to the argument in \S\,\ref{S11} (together with  the homomorphism (\ref{E4522})) that 
$\Phi_\alpha^*(\mcW_k^{2n})$ may be glued together to obtain an FC quantization 
\begin{align}
\mcW_\mcS
\end{align}
 on $(\mbA (\mcL_\mcS)^\times, \omega_\mcS)$.
 Consequently, we have obtained  the following assertion.
  (In the case $n=1$, one verifies that the symplectic structure $\omega_\mcS$ coincides with $\check{\omega}_\mbL$ and the asserted construction of FC quantizations  is consistent with  $\bigstar_{X, \mbL}$ mentioned in our main theorem.)

\vspace{3mm}
\bt  \label{T009h1gh}\leavevmode\\
 \ \ \ 
 Let $S$ be a smooth variety of dimension $2n$ (where $n$ is a positive integer).  
Then, by means of  a Frobenius-$\mr{Sp}$ structure  $\mcS^{\heartsuit \diamondsuit}$ on $S$, we can construct canonically a theta characteristic $\mbL := (\mcL_\mcS, \psi_\mcS)$ on $S$, a symplectic structure $\omega_\mcS$ on $\mbA (\mcL_\mcS)^\times$,  and an FC quantization $\mcW_\mcS$ on the resulting symplectic variety $(\mbA (\mcL_\mcS)^\times, \omega_\mcS)$.
 \et
\vspace{3mm}

\end{document}